\documentclass[STIX1COL]{WileyNJD-v2}
\usepackage{moreverb}
\usepackage{graphicx}
\usepackage{hyperref}  
\usepackage{amsmath}
\usepackage{amsfonts}
\usepackage{mathrsfs}
\usepackage{subfigure}
\usepackage{multirow}
\usepackage{mwe}
\usepackage{makecell}
\usepackage{caption}
\usepackage{color}

\usepackage{dirtytalk}

\newtheorem{exm}{Example}[section]
\newtheorem{lem}{Lemma}[section]

\newcommand\BibTeX{{\rmfamily B\kern-.05em \textsc{i\kern-.025em b}\kern-.08em
T\kern-.1667em\lower.7ex\hbox{E}\kern-.125emX}}

\articletype{Preprint}%

\begin{document}

\title{A new approach to solving multi-order fractional equations using BEM and Chebyshev matrix}

\author[1,2]{Moein Khalighi*}

\author[1,3]{Mohammad Amirianmatlob}


\author[1]{Alaeddin Malek}

\authormark{Moein Khalighi \textsc{et al}}

\address[1]{\orgdiv{Department of Applied Mathematics}, \orgname{Tarbiat Modares University}, \orgaddress{\state{Tehran}, \country{Iran}}}

\address[2]{\orgdiv{Department of Future Technologies}, \orgname{University of Turku}, \orgaddress{\state{Turku}, \country{Finland}}}

\address[3]{\orgdiv{Department of Mathematics and Statistics}, \orgname{Dalhousie University}, \orgaddress{\state{Halifax}, \country{Canada}}}


\corres{*Moein Khalighi (moein.khalighi@utu.fi)}

\presentaddress{Department of Future Technologies, Faculty of Sciences and Engineering, Yliopistonmaki, FI-20014, University of Turku, Finland}

\abstract[Abstract]{In this paper, the boundary element method is combined with Chebyshev operational matrix technique to solve two-dimensional multi-order time-fractional partial differential equations; nonlinear and linear in respect to spatial and temporal variables, respectively. Fractional derivatives are estimated by Caputo sense. Boundary element method is used to convert the main problem into a system of a multi-order fractional ordinary differential equation. Then, the produced system is approximated by Chebyshev operational matrix technique, ans its condition number is analyzed. Accuracy and efficiency of the proposed hybrid scheme are demonstrated by solving three different types of two-dimensional time fractional convection-diffusion equations numerically. The convergent rates are calculated for different meshing within the boundary element technique. Numerical results are given by graphs and tables for solutions and different type of error norms.}



\maketitle


\section{Introduction}\label{sec1}
Fractional order differential operators are the representative of non-local phenomena while many integer-order differential operators are mostly applied to examine local phenomena \cite{teodoro2019review}; Therefore, fractional calculus can be useful to describe many of real-world problems which cannot be covered in the classic mathematical literature \cite{metzler2000random, Matlob2017}. Since the next state of many systems depend on its current and historical states, there is a great demand to improve topical methods for the real life problems~\cite{sun2018new, almeida2016modeling}. These problems happen in bioengineering \cite{bioen}, solid mechanics \cite{solid}, anomalous transport \cite{metzler}, continuum and statistical mechanics \cite{mechanic}, economics \cite{eco}, relaxation  electrochemistry \cite{oldham}, diffusion procedures \cite{pod}, and complex networks \cite{saf,megh}, optimal control problems~\cite{opt3, opt2}. Fractional diffusion equations are largely used in describing abnormal convection phenomenon of liquid in medium.
Models of convection-diffusion quantities play significant roles in many practical applications~\cite{sun2018new, almeida2016modeling}, especially those involving fluid flow and heat transfer, such as thermal pollution in river system, leaching of salts in soils for computational simulations, oil reservoir simulations, transport of mass and energy, and global weather production. Numerical methods for convection-diffusion equations described by derivatives with integer order have been studied extensively \cite{malek2009}.
Due to the mathematical complexity, analytical solutions are very few and are restricted to the solution of simple fractional ordinary differential equations (ODEs) \cite{kilbas}. Several numerical techniques for solving fractional partial differential equations (PDEs) have been reported, such as variational iteration \cite{dehghan2011}, Adomian decomposition \cite{ford2009}, operational matrix of B-spline functions \cite{las2012}, operational matrix of Jacobi polynomials \cite{pod}, Jacobi collocation \cite{Lin}, operational matrix of Chebyshev polynomials \cite{malekCheby,bhra2013}, Legendre collocation \cite{mokhtar2015}, pseudo-spectral \cite{spect}, and operational matrix of Laguerre polynomials \cite{abd}, Pade approximation and two-sided Laplace transformations\cite{operationalmatrix}. Besides finite elements and finite differences \cite{Li}, spectral methods are one of the three main methodologies for solving fractional differential equations on computers \cite{rev}. The main idea of spectral methods is to express the solution of the differential equation as a sum of basis functions and then to choose the coefficients in order to minimize the error between the numerical and exact solutions as well as possible. Therefore, high accuracy and ease of implementing are two of the main features which have encouraged many researchers to apply such methods to solve various types of fractional integral and differential equation. In this article shifted chebyshev polynomials are used.\\
The boundary element method (BEM) always requires a fundamental solution to the original differential equation in order to avoid domain integrals in the formulation of the boundary integral equation. In addition, the nonhomogeneous and nonlinear terms are incorporated in the formulation by means of domain integrals. The basic idea of BEM is the transformation of the original differential equation into an equivalent integral equation only on the boundary, which has been widely applied in many areas of engineering such as fluid mechanics \cite{katsi2002}, magnetohydrodynamic \cite{malek2015}, and electrodynamics \cite{breb}.
\\In this paper the BEM is developed for the numerical solution of time fractional partial differential equations (TFPDEs) for nonhomogeneous bodies, which converts the main problem into a system of fractional ODE with initial conditions, described by an equation having a known fundamental solution. The proposed method introduces an additional unknown domain function, which represents the fictitious source function for the equivalent problem. This function is determined from a supplementary domain integral equation, which is converted to a boundary integral equation using a meshless technique based on global approximation by a radial basis function (RBF) series.
\\The Delaunay graph mapping method can be viewed as a fast interpolation scheme. Despite its efficiency, the mesh quality for large deformation may deteriorate near the boundary, in particular, if the deformation involves large rotation, which may even lead to an invalid Delaunay graph. Furthermore, the RBF method can generally better preserve the mesh quality near the boundary but the computational cost is much higher, as the mesh size increases. In order to develop methods that are more efficient and because of their flexibility and simplicity, the Delaunay graph based RBF method (DG-RBF) were proposed \cite{wang2015} to overcome the difficulty of meshing and remeshing the entire structure.
Thus, the pure boundary character of the method is maintained, since elements discretization and the integrations are limited only to the boundary. To obtain the fictitious source we use the Chebyshev spectral method based on operational matrix. The primary aim of this method is to propose a suitable way to approximate linear multi-order fractional ODEs with constant coefficients using a shifted Chebyshev Tau method, that guarantees an exponential convergence speed~\cite{taha}. Once the fictitious source is established, the solution of the original problem can be calculated from the integral representation of the solution in the substituted problem.\\
The outline of this paper is as follows. In Section~\ref{sec2}, we introduce the multi-order time fractional convection diffusion equation (TFCDE) for a class of the TFPDEs as a mathematical modelling of natural phenomena, and some basic preliminaries are also given. Section~\ref{sec3} is devoted to applying the BEM for converting the main problem into a system of multi-order fractional ODE with initial conditions. In Section~\ref{sec4}, the Chebyshev operational matrix (COM) of fractional derivative is obtained by applying the spectral methods to solve the generated multi-order fractional ODE. In order to demonstrate the efficiency and accuracy of the proposed method, along with the analysis of the condition number of COM, some numerical experiments are presented in Section~\ref{sec6} using the definitions and lemmas of Section~\ref{sec5}. Eventually, we conclude the paper with some remarks, and add the appendix including notation table~\ref{tabale: notation} to make more convenient understanding of the proposed algorithm.

\section{Problem statement}
\label{sec2}
 
Assume we are given the following initial boundary value problem for the multi-order time-fractional PDE in the two-dimensional domain $\Omega$ with boundary $\Gamma$,
\begin{equation}
\label{eq1}
\begin{gathered}
  \sum\limits_{j = 0}^k {{\gamma _j}({\mathbf{x}})D_c^{{\alpha _j}}u}  = A({\mathbf{x}}){u_{xx}} + 2B({\mathbf{x}}){u_{xy}} + C({\mathbf{x}}){u_{yy}} \hfill \\
  \quad \quad \quad \quad \quad \quad \quad + D({\mathbf{x}}){u_x} + E({\mathbf{x}}){u_y} + F({\mathbf{x}})u + g({\mathbf{x}},t), \hfill \\ 
\end{gathered} 
\end{equation}
where $A({\bf{x}})$, $B({\bf{x}})$, $C({\bf{x}})$, $D({\bf{x}})$, $E({\bf{x}})$, $F({\bf{x}})$ and ${\gamma _j}({\bf{x}})$ for $j=0,1,2,...,k$ and $g({\bf{x}},t)$ are specified functions their physical meaning depends on that of the field function $u({\bf{x}},t)$, and
\[{0 < {\alpha _0} < {\alpha _1} <  \cdots  < {\alpha _{k - 1}} < {\alpha _k},\quad\quad\,m - 1 < {\alpha _k} \leqslant m},\]
\[{{\mathbf{x}}=(x,y) \in \Omega \, \cup \,\Gamma ,\quad\quad\,t > 0},\]
subject to the boundary conditions
\begin{equation}
\label{eq2}
{\rm B}(u) = h({\mathbf{x}},t),\quad \quad {\mathbf{x}} \in \Gamma ,
\end{equation}
and the initial conditions
\begin{equation}
\label{eq3}
D_c^iu({\mathbf{x}},0) = {d_i}({\mathbf{x}}),\quad \quad i = 0,1,...,m - 1.
\end{equation}
In which $ m $ is an integer number and $ D_c^{{\alpha_j}}u $ is the Caputo fractional time derivative of order ${\alpha_j}$. The Caputo derivative \cite{pod}, is employed because initial conditions having direct physical meaning can be prescribed. This derivative is defined as 
\begin{equation}
\label{eq4}
D_c^\alpha u(x,y,t) = \left\{ \begin{gathered}
  \frac{1}{{\Gamma (m - \alpha )}}\int_0^t {\frac{{{u^{(m)}}(x,y,\tau ){\rm{d}\tau}}}{{{{(t - \tau )}^{1 + \alpha  - m}}}}} ,\quad m - 1 < \alpha  < m, \hfill \\
  \frac{{{d^m}}}{{d{t^m}}}u(x,y,t),\quad \quad \quad \quad \quad \quad \quad \alpha  = m \in \mathbb{N}. \hfill \\ 
\end{gathered}  \right.
\end{equation}
\\
${\rm B}( \cdot )$ is a linear operator with respect to spatial variables $x$, $y$ of order one. $h({\bf{x}},t)$ and ${d_i}({\bf{x}})$ $(i = 0,...,m - 1)$ are specified functions in Eq. \eqref{eq2} and \eqref{eq3}, respectively. It seems that we could be able to recover the multi-term of classical diffusion equation for ${\alpha _k} = 1$, ${\gamma _j}({\bf{x}}) = 0$ for $j = 0,...,k - 1$, ${\gamma _k}({\bf{x}}) \ne 0$  and the classical wave equation in presence of viscous damping for ${\alpha _k} = 2$, ${\alpha _{k - 1}} = 1$, ${{\gamma _k}({\bf{x}}) \ne 0}$, ${\gamma _{k - 1}}({\bf{x}}) \ne 0$ and ${\gamma _j}({\bf{x}}) = 0$ for $j = 0,...,k - 2$.
\section{Implementation of boundary element method}
\label{sec3}
Taking advantage of the following boundary element, the initial boundary value of the equation (\ref{eq1}-\ref{eq3}) is reformed into an ODE problem.
\\
Let $ u({\bf{x}},t)$ be the sought solution to the problem (\ref{eq1}-\ref{eq3}) and assume that $u$ is twice continuously differentiable in $\Omega$. After applying Laplace operator on $u$ we have \cite{katsi2002}
\begin{equation}
\label{eq5}
{\nabla ^2}u({\bf{x}},t) = \mathfrak{B}({\bf{x}},t),
\end{equation}
where $\mathfrak{B}({\bf{x}},t)$ known as an unknown fictitious source function. That is the solution of Eq. \eqref{eq1} could be established by solving Eq. \eqref{eq5} under the boundary condition \eqref{eq2}, if convenient $\mathfrak{B}({\bf{x}},t)$ is first established. This is accomplished by adhering to the following procedure. \\

We write the solution of Eq. \eqref{eq5} in the integral form. Thus, for $u_n$ as the normal derivative of $u$ we have \cite{katsi2011}
\begin{equation}
\label{eq6}
\varepsilon u({\bf{x}},t) = \int_\Omega  {{u^*}\mathfrak{B} \mathrm{d}\Omega  - } \int_\Gamma  {({u^*}{u_n} - u_n^*u)\mathrm{d}\Gamma },  \quad \quad {\bf{x}} \in \Omega  \cup \Gamma, 
\end{equation}
where ${u^*} = \ln r/2\pi $ is the fundamental solution to Eq. \eqref{eq5}, $r$ is the distance between any two points and also $u_n^*$ stands for its  normal derivative on the boundary. $\varepsilon $ is the free term coefficient taking the values $\varepsilon=1$ if ${\bf{x}} \in \Omega$, $\varepsilon=\theta/{2\pi}$ if ${\bf{x}} \in \Gamma$, otherwise $\varepsilon=0$; $\theta$ is the interior angle between the tangents of boundary at point $\bf{x}$. $\varepsilon=1/2$ for points where the boundary is smooth. After applying 
Eq. \eqref{eq6}, to boundary points by means of Greens second identity \cite{katsi2005}, we yield the boundary integral equation
\begin{equation}
\label{eq7}
\frac{\theta }{{2\pi }}u({\bf{x}},t) = \sum\limits_{j = 1}^M {{b_j}\left[ {\frac{1}{2}{{\hat u}_j} + \int_\Gamma  {\left( {{u^*}{{({{\hat u}_n})}_j} - u_n^*{{\hat u}_j}} \right)\mathrm{d}\Gamma } } \right] - } \int_\Gamma  {({u^*}{u_n} - u_n^*u)\mathrm{d}\Gamma }, \; {\bf{x}} \in \Gamma ,
\end{equation}
where ${\hat u_j}$ is a certain solution of the equation
\begin{equation}
    {\nabla ^2}{\hat u_j} = {f_j},\quad \quad j = 1,2,...,M,
\end{equation}
Also $M$ is the number of interior points inside $\Omega$. Here $b_j$ are the coefficients that must be determined to satisfy
\[\mathfrak{B}({\bf{x}},t) = \sum\limits_{j = 1}^M {{b_j}{f_j}}, \]
where ${f_j} = {f_j}(r)$, $r = \left| {{\bf{x}} - {{\bf{x}}_j}} \right|$ is a set of radial basis approximating functions; ${{\bf{x}}_j}$ are collocation points in $\Omega$. The radial basis function method is used to map the nodes rather than that based on surface or volume ratios \cite{wang2015}. The algorithm is set out in the following procedure; At first we generate the Delaunay graph by using all the boundary nodes of the original mesh, and then we locate the mesh points in the graph, after that we move the Delaunay graph according to the specified geometric motion/deformation, and the final step is mapping the mesh points in the new graph according to the RBF matrix and Delaunay triangle. Procedures before the last step are exactly the same as the original Delaunay graph mapping method \cite{dela}; hence the details of these steps are skipped in this paper. The difference is in the last step, where the radial basis function interpolation is used to calculate the displacement of the internal mesh nodes from the given displacement of the Delaunay triangle nodes on the boundary, while the original Delaunay mapping method uses surface or volume ratios to calculate the displacement of inner nodes.
Eq. \eqref{eq7} is solved numerically by using the BEM. The boundary integrals in Eq. \eqref{eq7} are approximated using $N$ boundary nodal points. Here $\varepsilon=1/2$, as we ace the smooth boundary. 

The domain integral can be evaluated when the fictitious source is estimated by a radial basis function series and subsequently it is reformed to a boundary line integral using the Green’s reciprocal identity \cite{katsi2005}. For the sake of simplification, we use multiquadric radial basis function in practice. $M$ internal nodals are here used to define Delaunay linear triangular elements in $\Omega$. Therefore, after discretization and application of the boundary integral equation \eqref{eq7} at the $N$ boundary nodal points we have
\begin{equation}
\label{eq8}
{\bf{Hu}} = {\bf{G}}{{\bf{u}}_n} + {\bf{{\mathbf{\bar A}}b}},
\end{equation}
where $\bf{H}$ and $\bf{G}$ are $N \times N$ known as coefficient matrices originating from the integration of the kernel functions on the boundary elements and ${\mathbf{\bar A}}$ is an $N \times M$ coefficient matrix originating from the integration of the kernel function on the domain elements; $\bf{u}$ and ${\bf{u}}_n$ are vectors containing the nodal values of the boundary displacements and their normal derivatives. Also, $\bf{b}$ is the vector of the nodal values of the fictitious source at the $M$ internal nodal points. \\
For a second order differential equation, the boundary condition is a correlation of  ${\delta _1}({\bf{x}})u + {\delta _2}({\bf{x}}){u_n} = h({\bf{x}},t)$; after applying it at the $N$ boundary nodal points yields
\begin{equation}
\label{eq9}
{{\pmb{\delta }}_1}{\bf{u}} + {{\pmb{\delta }}_2}{{\bf{u}}_n} = {\bf{h}}(t),
\end{equation}
where ${\pmb{\delta}}_1$ and ${\pmb{\delta}}_2$ are $N \times N$ known diagonal matrices and
 ${\bf{h}}(t)=$
 $[ {h\left( {{x_{B{P_1}}},t} \right), \ldots ,}$
 $h({x_{B{P_N}}},t){]^T}$
 is a known boundary vector, where ${x_{B{P_j}}}$ are $N$ boundary nodal points. Eqs. \eqref{eq8} and \eqref{eq9} can be combined and solved for $\bf{u}$ and ${\bf{u}}_n$. This yields
\begin{equation}	
\label{eq10}
\begin{array}{l}
{\bf{u}} = \left[ {{{({{\pmb{\delta }}_{{1}}}{\bf{ + }}{{\pmb{\delta }}_{{2}}}{{\bf{G}}^{{{ - 1}}}}{\bf{H}})}^{{{ - 1}}}}({{\pmb{\delta }}_{{2}}}{{\bf{G}}^{{{ - 1}}}}{\mathbf{\bar A}})} \right]{\bf{b}} + {({{\pmb{\delta }}_{{1}}} + {{\pmb{\delta }}_{{2}}}{{\bf{G}}^{{{ - 1}}}}{\bf{H}})^{{{ - 1}}}}{\bf{h}}(t),\\
{{\bf{u}}_n} = \left[ {{{({{\pmb{\delta }}_1}{{\bf{H}}^{{{ - 1}}}}{\bf{G}} + {{\pmb{\delta }}_{{2}}})}^{{{ - 1}}}}( - {{\pmb{\delta }}_1}{{\bf{H}}^{{{ - 1}}}}{\mathbf{\bar A}})} \right]{\bf{b}} + {({{\pmb{\delta }}_1}{{\bf{H}}^{{{ - 1}}}}{\bf{G}} + {{\pmb{\delta }}_{{2}}})^{{{ - 1}}}}{\bf{h}}(t),
\end{array}
\end{equation}
Further, differentiating Eq. \eqref{eq6} for points inside the domain ($\varepsilon=1$) with respect to $x$ and $y$, using the same discretization and collocating at the $M$ internal nodal points, we have the following expression for the spatial derivatives
\begin{equation}
\label{eq11}
{{\bf{\hat u}}_{pq}} = {{\bf{\hat H}}_{pq}}{\bf{u}} + {{\bf{\hat G}}_{pq}}{{\bf{u}}_n} + {{\bf{\hat A}}_{pq}}{\bf{b}},\quad\quad\quad p,q = 0,x,y
\end{equation}
where the ${{\bf{\hat u}}_{pq}}$ is vector of values for $u$ and its derivatives at the $M$ internal nodal points; ${\bf{\hat H}}_{pq}$ and ${\bf{\hat G}}_{pq}$ are $M \times N$ known coefficient matrices originating from the integration of the kernel functions on the boundary elements and ${\bf{\hat A}}_{pq}$ is an $M \times M$ coefficient matrix originating from the integration of the kernel functions on the domain elements.\\
Eliminating  $\bf{u}$ and ${\bf{u}}_n$ from Eq. \eqref{eq11} using Eqs. \eqref{eq10} yields 
\begin{equation}
\label{eq12}
{{\bf{\hat u}}_{pq}} = {{\bf{U}}_{pq}}{\bf{b}} + {{\bf{c}}_{pq}},\quad \quad p,q = 0,x,y
\end{equation}
where

\begin{equation}\label{eq:Uc}
\begin{array}{l}
{{\bf{U}}_{pq}} = {{{\bf{\hat H}}}_{pq}}{({{\pmb{\delta }}_{{1}}}{{ + }}{{\pmb{\delta }}_{{2}}}{{\bf{G}}^{{{ - 1}}}}{\bf{H}})^{{{ - 1}}}}({{\pmb{\delta }}_{{2}}}{{\bf{G}}^{{{ - 1}}}}{\mathbf{\bar A}}) + {{{\bf{\hat G}}}_{pq}}{({{\pmb{\delta }}_1}{{\bf{H}}^{{{ - 1}}}}{\bf{G}} + {{\pmb{\delta }}_{{2}}})^{{{ - 1}}}}( - {{\pmb{\delta }}_1}{{\bf{H}}^{{{ - 1}}}}{\mathbf{\bar A}}) + {{{\bf{\hat A}}}_{pq}},\\
{{\bf{c}}_{pq}} = \left[ {{{{\bf{\hat H}}}_{pq}}{{({{\pmb{\delta }}_{{1}}}{{ + }}{{\pmb{\delta }}_{{2}}}{{\bf{G}}^{{{ - 1}}}}{\bf{H}})}^{{{ - 1}}}} + {{{\bf{\hat G}}}_{pq}}{{({{\pmb{\delta }}_1}{{\bf{H}}^{{{ - 1}}}}{\bf{G}} + {{\pmb{\delta }}_{{2}}})}^{{{ - 1}}}}} \right]{\bf{h}}(t),
\end{array}
\end{equation}

The final step of the method is to apply Eq. \eqref{eq1} at the $M$ internal nodal points. This gives 
\begin{equation}
\label{eq13}
\sum\limits_{j = 0}^k {{{\pmb{\gamma }}_j}D_c^{{\alpha _j}}{\bf{\hat u}}}  = {\bf{A}}{{\bf{\hat u}}_{xx}} + 2{\bf{B}}{{\bf{\hat u}}_{xy}} + {\bf{C}}{{\bf{\hat u}}_{yy}} + {\bf{D}}{{\bf{\hat u}}_x} + {\bf{E}}{{\bf{\hat u}}_y} + {\bf{F\hat u}} + {\bf{g}}(t),
\end{equation}
where $\hat{\bf{u}}={\hat{\bf{u}}}_{00}$ and ${\pmb{\gamma}}_j$, $\bf{A}$, $\bf{B}$, $\bf{C}$, $\bf{D}$, $\bf{E}$ and $\bf{F}$ are $M \times M$ known diagonal matrices including the nodal values of the corresponding functions $\gamma_{j}(\bf{x})$ $A(\bf{x})$, $B(\bf{x})$, $C(\bf{x})$, $D(\bf{x})$, $E(\bf{x})$ and $F(\bf{x})$, respectively, and ${\bf{g}}(t)={\left[ {g\left( {{x_{I{P_1}}},t} \right), \ldots ,g\left( {{x_{I{P_M}}},t} \right)}\right]^T}$ is a known internal vector, where ${x_{I{P_j}}}$ are $M$ internal nodal points. Substituting the corresponding terms from Eq. \eqref{eq12} into Eq. \eqref{eq13} yields
\begin{equation}
\label{eq14}
\sum\limits_{j = 0}^k {{{\bf{S}}_j}D_c^{{\alpha _j}}{\bf{b}}(t)}  = {\bf{Nb}}(t) + {\bf{f}}(t),
\end{equation}
where 
\begin{equation}\label{eq:SNF}
\begin{array}{l}
{{\bf{S}}_j} = {{\pmb{\gamma }}_j}{\bf{U}},\\
{\bf{N}} = {\bf{A}}{{\bf{U}}_{xx}} + 2{\bf{B}}{{\bf{U}}_{xy}} + {\bf{C}}{{\bf{U}}_{yy}} + {\bf{D}}{{\bf{U}}_x} + {\bf{E}}{{\bf{U}}_y} + {\bf{FU}},\\
{\bf{f}}(t) = {\bf{A}}{{\bf{c}}_{xx}} + 2{\bf{B}}{{\bf{c}}_{xy}} + {\bf{C}}{{\bf{c}}_{yy}} + {\bf{D}}{{\bf{c}}_x} + {\bf{E}}{{\bf{c}}_y} + {\bf{Fc}} + {\bf{g}}(t) - \mathop \sum \limits_{j = 0}^k {{\pmb{\gamma }}_j}D_c^{{\alpha _j}}({\bf{c}}),
\end{array}
\end{equation}
in which $\bf{U}=\bf{U_{00}}$ and $\bf{c}=\bf{c_{00}}$ for $j = 0,...,k$. Now, from Eq. \eqref{eq12}, we can write the initial conditions \eqref{eq3} for $i = 0,1,...,m - 1$ in the form
\begin{equation}
\label{eq16}
{{\bf{b}}^{(i)}}(0) = {{\bf{U}}^{ - 1}}\left[ {{{\bf{d}}_i} - {{\bf{c}}^{(i)}}(0)} \right],
\end{equation}
where ${{\bf{d}}_i}(t)={\left[ {{{d}_i}\left( {{x_{I{P_1}}},t} \right), \ldots ,{{d}_i}\left( {{x_{I{P_M}}},t} \right)}\right]^T}$.\\
The above proposed procedure reduces the problem of multi-order two-dimensional time fractional PDE (\ref{eq1}-\ref{eq3}) to a simpler system of multi-term fractional ODE \eqref{eq14} with initial condition \eqref{eq16}. The existence, uniqueness, and continuous dependence of the system~\eqref{eq14}-\eqref{eq16} when $\textbf{S}_{k}=1$ can be rigorously discussed (see e.g. Diethelm and Neville's paper~\cite{DIETHELM2004621}). In the next section, we show the implementation of Chebyshev operational matrix, as a spectral technique~\cite{rev} for fractional calculus, to solve the system of initial value problem \eqref{eq14}-\eqref{eq16}.

\section{COM method for system of multi-order fractional ODEs}
\label{sec4}
The Chebyshev polynomials $T_{i}(t)$ are defined on the interval $(-1,1)$~\cite{taha}. Thus, by changing variable $t\rightarrow\frac{2t}{L}-1$, the shifted Chebyshev polynomials ${T_{L,i}}(t)$ of degree $i$ on the interval $t \in (0,L)$, with an orthogonality relation can be introduced by~\cite{rev, DOHA20112364}
\[{T_{L,i}}(t) = i\sum\limits_{j = 0}^i {{{( - 1)}^{i - j}}\frac{{(i + j - 1)!{2^{2j}}}}{{(i - j)!(2j)!{L^j}}}{t^j}}, \]
where ${T_{L,i}}(0) = {( - 1)^i}$ and ${T_{L,i}}(L) = 1$. 
In this form, ${T_{L,i}}(t)$ may be generated by the following recurrence formula: 
\begin{equation}\label{eq: collocation}
{T_{L,i + 1}}(t) = 2({2t}/{L} - 1){T_{L,i}}(t) - {T_{L,i - 1}}(t),\quad \quad \quad i = 1,2,...
\end{equation}
where ${T_{L,0}}(t) = 1$ and ${T_{L,1}}(t) = {2t}/{L} - 1$. Therefore, a given function $f\in{L^{2}[0,1]}$ may be approximated by $K+1$ terms of shifted Chebyshev polynomials as \[f(t)\simeq{f_{K}(t)}=\sum_{i=0}^{K}{c_{i}T_{L,i}(t)},\] where the coefficients $c_{i}$ are described by weight functions $w_{L}(t)=\frac{1}{\sqrt{Lt-t^2}}$ as $c_{j}=\frac{1}{h_{i}}\int_{0}^{L}f(t)T_{L,i}(t)w_{L}(t)dt$; in which $h_{i}=\pi$ for $i=0$, otherwise $h_{i}=\frac{\pi}{2}$. If we set
\begin{equation}\label{eq:phi}
\Phi (t) = {\left[ {{T_{L,0}}(t),{T_{L,1}}(t),...,{T_{L,K}}(t)} \right]^T},
\end{equation}
and suppose $\upsilon  > 0$ and the ceiling function $\left\lceil \upsilon  \right\rceil $ denotes the smallest integer greater than or equal to $\upsilon$, then
\begin{equation}
\label{eq27}
D_c^\upsilon \Phi (t)\simeq {\mathfrak{D}^{(\upsilon )}}\Phi (t),
\end{equation}
where ${\mathfrak{D}^{(\upsilon )}}$ is the $(K+1)\times(K+1)$ COM of derivatives of order $\upsilon$ in the Caputo sense and is defined by~\cite{rev, DOHA20112364}: 
\begin{equation}
\label{eq28}
{\mathfrak{D}}^{(\upsilon )} = \left( {\begin{array}{*{20}{c}}
0&0&0& \cdots &0\\
 \vdots & \vdots & \vdots & \cdots & \vdots \\
0&0&0& \cdots &0\\
{{S_\upsilon }\left( {\left\lceil \upsilon  \right\rceil ,0} \right)}&{{S_\upsilon }\left( {\left\lceil \upsilon  \right\rceil ,1} \right)}&{{S_\upsilon }\left( {\left\lceil \upsilon  \right\rceil ,2} \right)}& \cdots &{{S_\upsilon }\left( {\left\lceil \upsilon  \right\rceil ,K} \right)}\\
 \vdots & \vdots & \vdots & \cdots & \vdots \\
{{S_\upsilon }\left( {i,0} \right)}&{{S_\upsilon }\left( {i,1} \right)}&{{S_\upsilon }\left( {i,2} \right)}& \cdots &{{S_\upsilon }\left( {i,K} \right)}\\
 \vdots & \vdots & \vdots & \cdots & \vdots \\
{{S_\upsilon }\left( {K,0} \right)}&{{S_\upsilon }\left( {K,1} \right)}&{{S_\upsilon }\left( {K,2} \right)}& \cdots &{{S_\upsilon }\left( {K,K} \right)}
\end{array}} \right)
\end{equation}
where
\[{S_\upsilon }(i,j) = \sum\limits_{\lambda = \left\lceil \upsilon  \right\rceil }^i {\frac{{{{( - 1)}^{i - \lambda}}2i(i + \lambda - 1)!\Gamma \left( {\lambda - \upsilon  + \frac{1}{2}} \right)}}{{{\rho _j}{L^\upsilon }\Gamma \left( {\lambda + \frac{1}{2}} \right)\left( {i - \lambda} \right)!\Gamma (\lambda - \upsilon  - j + 1)\Gamma (\lambda + j - \upsilon  + 1)}}} ,\]
where ${\rho _0} = 2$, ${\rho _\lambda} = 1$, $\lambda \ge 1$. Note that in ${\mathfrak{D}}^{(\upsilon )}$, the first $\left\lceil \upsilon  \right\rceil $ rows are all zero.
In order to solve Eq. \eqref{eq14} with initial conditions \eqref{eq16}, we approximate ${\bf{b}}(t)$ and ${\bf{f}}(t)$ in terms of shifted Chebyshev polynomials as
\begin{equation}
\label{eq29}
{\bf{b}}(t) = \left[ {\begin{array}{*{20}{c}}
{\sum\limits_{i = 0}^K {\psi _i^1{T_{L,i}}(t)} }\\
 \vdots \\
{\sum\limits_{i = 0}^K {\psi _i^M{T_{L,i}}(t)} }
\end{array}} \right] = \left[ {\begin{array}{*{20}{c}}
{{\Psi ^1}^{^T}\Phi (t)}\\
 \vdots \\
{{\Psi ^M}^{^T}\Phi (t)}
\end{array}} \right] = {\pmb{\Psi }}\Phi (t),
\end{equation}
\begin{equation}
\label{eq30}
{\bf{f}}(t) = \left[ {\begin{array}{*{20}{c}}
{\sum\limits_{i = 0}^K {\ell _i^1{T_{L,i}}(t)} }\\
 \vdots \\
{\sum\limits_{i = 0}^K {\ell _i^M{T_{L,i}}(t)} }
\end{array}} \right] = \left[ {\begin{array}{*{20}{c}}
{{\omega ^1}^{^T}\Phi (t)}\\
 \vdots \\
{{\omega ^M}^{^T}\Phi (t)}
\end{array}} \right] = {\pmb{\omega }}\Phi (t),
\end{equation}
where for $j=1,...,M$
\[\begin{array}{l}
{\Psi ^j}^{^T} = \left[ {\psi _0^j,...,\psi _K^j} \right],\\
{\omega ^j}^{^T} = \left[ {\ell _0^j,...,\ell _K^j} \right],
\end{array}\]
and $\pmb{\Psi}$, $\pmb{\omega}$ are $M\times(K+1)$ matrices that are defined as 
\[{\pmb{\Psi }} = \left[ {\begin{array}{*{20}{c}}
{{\Psi ^1}^{^T}}\\
 \vdots \\
{{\Psi ^M}^{^T}}
\end{array}} \right],\quad \quad \quad \quad {\pmb{\omega }} = \left[ {\begin{array}{*{20}{c}}
{{\omega ^1}^{^T}}\\
 \vdots \\
{{\omega ^M}^{^T}}
\end{array}} \right].\]
For $j=0,...,k$, Eq. \eqref{eq27} and Eq. \eqref{eq29} can be used to write
\begin{equation}
\label{eq31}
D_c^{{\alpha _j}}{\bf{b}}(t)\simeq {\pmb{\Psi }}D_c^{{\alpha _j}}\Phi (t)\simeq {\pmb{\Psi }}{\mathfrak{D}^{({\alpha _j})}}\Phi (t).
\end{equation}
Employing Eqs.\eqref{eq29}, \eqref{eq30} and \eqref{eq31}, then the residual for Eq. \eqref{eq14} can be written as 
\begin{equation}
\label{eq32}
{\bf{R}}(t) = \left( {\sum\limits_{j = 0}^k {{{\bf{S}}_j}{\pmb{\Psi }}} {\mathfrak{D}^{({\alpha _j})}} - {\bf{N}\pmb{\Psi }} - {\pmb{\omega }}} \right)\Phi (t).
\end{equation}
${\bf{R}}(t)$ is a $M$ vector with respect to $t$. If $R_j(t)$ be the $j$th component of ${\bf{R}}(t)$, then in a typical Tau method \cite{taha}, we generate $M(K-m+1)$ linear equations with $M(K+1)$ unknown coefficients of	$\pmb{\Psi}$ by applying 
\begin{equation}
\label{eq33}
\left\langle {{R_i}(t),{T_{L,j}}(t)} \right\rangle  = \int_0^L {{R_i}(t){T_{L,j}}(t)\mathrm{d}t = 0} ,\quad i = 1,...,M,\quad j = 0,1,...,K - m.
\end{equation}
Also, by substituting Eq. \eqref{eq29} into Eq. \eqref{eq16}, and with the fact that 
\[{\mathfrak{D}^{(n)}} = {({\mathfrak{D}^{(1)}})^n},\quad n \in \mathbb{N},\]
we get a system of $Mm$ linear equations with $M(K+1)$ unknown coefficients for $\pmb{\Psi}$ as following 
\begin{equation}
\label{eq34}
{{\mathbf{b}}^{(i)}}(0) = {\pmb{\Psi }}{\mathfrak{D}^{(i)}}\Phi (0) = {{\mathbf{U}}^{ - 1}}\left[ {{{\mathbf{d}}_i} - {{\mathbf{c}}^{(i)}}(0)} \right],\quad \quad i = 0,1,...,m - 1.
\end{equation}
Equations \eqref{eq33} and \eqref{eq34} can be rewritten in the matrix form
\begin{equation}\label{matrixform}
    \pmb{A}\pmb{\Psi}=\pmb{B},
\end{equation}
where $\pmb{A}$ is an $M(K+1)\times M(K+1)$ coefficient matrix. The system of algebraic equations \eqref{matrixform} can be easily solved for the unknown vector $\pmb{\Psi}$. Consequently, ${\bf{b}}(t)$ given in Eq. \eqref{eq29} can be calculated, which gives a solution of Eq. \eqref{eq14} with the initial conditions \eqref{eq16}. Once the vector ${\bf{b}}(t)$ of the values of the fictitious source at the $M$ internal nodal points has been established, then the solution of Eq. \eqref{eq1} and its derivatives can be computed from Eq. \eqref{eq12}. For the points ${\bf{x}}=(x,y)$ that do not coincide with the prespecified internal nodal points, the solution could be drawn from the discretized counterpart of Eq. \eqref{eq6} with $\varepsilon=1$ using the same boundary and new domain discretization. Note that here the matrices ${{\mathbf{\hat H}}_{pq}}$, ${{\mathbf{\hat G}}_{pq}}$ and ${{\mathbf{\hat A}}_{pq}}$ corresponding to previous internal nodes plus the additional points must be recomputed.

\section{Error estimation}\label{sec5}
 The convergence of the proposed method is shown by employing the following error norms, maximum error (${L_\infty} $), maximum relative error ($MRE$) to assess the accuracy of the method in multi-scale problems, and root  mean square ($RMS$) to  globally examine the method efficiency,
\begin{equation}
\label{Linf}
{L_\infty } = \mathop {\max }\limits_{1 \leqslant i \leqslant M} \left| {u_{ex}^i - u_{app}^i} \right|,
\end{equation}
\begin{equation}
\label{MRE}
MRE = \mathop {\max }\limits_{1 \leqslant i \leqslant M} \left| {\frac{{u_{ex}^i - u_{app}^i}}{{u_{ex}^i}}} \right|,
\end{equation}
\begin{equation}
\label{RMS}
 RMS = \sqrt {\frac{1}{M}\sum\limits_{i = 1}^M {{{(u_{ex}^i - u_{app}^i)}^2}} } ,
\end{equation}
where $u_{ex}^i$ and $u_{app}^i$ denote the $i$th components of the exact and approximated solutions, respectively, and $M$ denotes the number of internal points. It is not convenient to certainly determine what is the convergence rate of the proposed hybrid method; for example, for the number of nodal points $N$ and $M$, and the size of the shifted Chebyshev polynomials $K$, the convergence order of the method would be $O(f({r^p},{{\tau}^q}))$ where $f$ is a function of the convergence rate of BEM with the order $p$ and the convergence rate of COM with the order $q$. The accuracy of the method depends on several factors, the convergence speed for BEM, the domain and boundary discretization, the shape parameters of radial basis functions, the orders of the derivatives, and condition number of COMs, as such the analysis of truncation errors for methods solving a two dimensional multi-term time fractional differential equation is not straightforward. Nonetheless, information about matrix $\pmb{A}$ from the algebraic system \eqref{matrixform}, particularly its condition number, will be useful. The condition number is defined by\cite{taha}
\[Cond(\pmb{A})=\frac{\mathrm{max}\{|\mathbf{\lambda}|:\mathrm{det}(\pmb{A}-\lambda I)=0 \}}{\mathrm{min}\{|\mathbf{\lambda}|:\mathrm{det}(\pmb{A}-\lambda I)=0 \}},\]
such that a matrix with a large condition number is so-called \textit{ill conditioned}, whereas the matrix is named \textit{well conditioned} if its condition number is of a moderate size. We also suggest two $P$-orders in the following lemmas to examine the rate of convergence for BEM and COM distinctly. The first one is directly tested by the exact solution and the effect of domain-discreitization. While the second one is addressed by comparing a sequence of numerical solutions of the ODE system~\eqref{eq14} with different degree sizes of COMs which have been offered exponential rates of convergence accuracy for smooth problems in simple geometries~\cite{taha}.

\begin{lem}\label{lemma1}
Let the vector ${\mathbf{U}}$ be the exact solution of the initial boundary value problem (\ref{eq1}-\ref{eq3}) and ${{\mathbf{u}}_1}$,${{\mathbf{u}}_2}$ the approximate solutions with $N_1$, $M_1$ and $N_2$, $M_2$ of nodal points, respectively. Then, the computational order of the BEM method proposed in Section \ref{sec3} can be calculated with
${P_{r}}$-order$\simeq{{\log \left( {{{{E_{r_1}}} \mathord{\left/
 {\vphantom {{{E_{r_1}}} {{E_{r_2}}}}} \right.
 \kern-\nulldelimiterspace} {{E_{r_2}}}}} \right)} \mathord{\left/
 {\vphantom {{\log \left( {{{{E_{r_1}}} \mathord{\left/
 {\vphantom {{{E_{r_1}}} {{E_{r_2}}}}} \right.
 \kern-\nulldelimiterspace} {{E_{r_2}}}}} \right)} {\log \left( {{{{{r} _1}} \mathord{\left/
 {\vphantom {{{{r}_1}} {{{r}_2}}}} \right.
 \kern-\nulldelimiterspace} {{{r}_2}}}} \right)}}} \right.
 \kern-\nulldelimiterspace} {\log \left( {{{{{r}_1}} \mathord{\left/
 {\vphantom {{{{r}_1}} {{{r}_2}}}} \right.
 \kern-\nulldelimiterspace} {{{r}_2}}}} \right)}}$
in which $E_{r_1}$ and $E_{r_2}$ are corresponded $RMS$ errors \eqref{RMS} with the relative boundary mesh size ${r}_1=1/{N_1}$ and ${r}_2=1/{N_2}$, respectively.
\end{lem}
\begin{proof}
When the leading terms in the spatial-discretization error are proportional to ${r_1}^p$ and ${r_2}^p$, and ${\left\| . \right\|_{RMS}}$ denoting the root mean square norm \eqref{RMS},
\[{\left\| {{\mathbf{U}} - {{\mathbf{u}}_1}} \right\|_{RMS}} = {E_{{r}_1}} = c_{1}{r}_1^p\simeq{c{r}_1^p},\quad {\left\| {{\mathbf{U}} - {{\mathbf{u}}_2}} \right\|_{RMS}} = {E_{{r}_2}} = c_2{r}_2^p\simeq{c{r}_2^p}.\]
Hence
\[\frac{{{E_{{r}_1}}}}{{{E_{{r}_2}}}} \simeq \frac{{{r_1}^p}}{{{r_2}^p}},\]
then taking logarithm from both sides yields
\begin{equation*}
p\simeq\frac{{\log \left( {\frac{{{E_{{r}_1}}}}{{{E_{{r}_2}}}}} \right)}}{{\log \left( {\frac{{{{r}_1}}}{{{{r}_2}}}} \right)}}. \qed
\end{equation*}

\end{proof}

\begin{lem}\label{lemma2}
Let the vector ${{\mathbf{b}}_{ex}}$ be the exact solution and ${{\mathbf{b}}_1},{{\mathbf{b}}_2}$ and ${{\mathbf{b}}_3}$ be the approximate solutions of the multi-term fractional ODE \eqref{eq14} with the initial condition \eqref{eq16} at the same $M$ fictitious source points using $K_1$, $K_2$, and $K_3$ the numbers of shifted Chebyshev polynomials, respectively. With considering this proportion 
\begin{equation}
\label{propo}
\frac{{{K_1}}}{{{K_2}}} = \frac{{{K_2}}}{{{K_3}}}, 
\end{equation}
the temporal convergence order for the COM method presented in Section \ref{sec4} is estimated using 
${P_\tau }$-order$\simeq{{\log \left( {\frac{{{{\left\| {{{\mathbf{b}}_2} - {{\mathbf{b}}_1}} \right\|}_b}}}{{{{\left\| {{{\mathbf{b}}_3} - {{\mathbf{b}}_2}} \right\|}_b}}}} \right)} \mathord{\left/
 {\vphantom {{\log \left( {\frac{{{{\left\| {{{\mathbf{b}}_2} - {{\mathbf{b}}_1}} \right\|}_b}}}{{{{\left\| {{{\mathbf{b}}_3} - {{\mathbf{b}}_2}} \right\|}_b}}}} \right)} {\log \left( {\frac{{{\tau _1}}}{{{\tau _2}}}} \right)}}} \right.
 \kern-\nulldelimiterspace} {\log \left( {\frac{{{\tau _1}}}{{{\tau _2}}}} \right)}}$
in which the norm
${\left\| . \right\|_b}$
define as
${\left\| {{{\mathbf{b}}_s} - {{\mathbf{b}}_t}} \right\|_b} = \sqrt {\frac{1}{M}\sum\limits_{i = 1}^M {{{\left( {b_s^i - b_t^i} \right)}^2}} } $
where $b_s^i$ and $b_t^i$ denote the $i$th components, and 
${\tau _1} = {1 \mathord{\left/
 {\vphantom {1 {{K_1}}}} \right.
 \kern-\nulldelimiterspace} {{K_1}}}$,
 ${\tau _2} = {1 \mathord{\left/
 {\vphantom {1 {{K_2}}}} \right.
 \kern-\nulldelimiterspace} {{K_2}}}$,
 and 
 ${\tau _3} = {1 \mathord{\left/
 {\vphantom {1 {{K_3}}}} \right.
 \kern-\nulldelimiterspace} {{K_3}}}$.
 \end{lem}
 \begin{proof} When the leading terms in the error of COM are proportional to ${\tau_1}^p$, ${\tau_2}^p$ and ${\tau_3}^p$, 
\[{\left\| {{{\mathbf{b}}_{ex}} - {{\mathbf{b}}_1}} \right\|_b} = c_{1}'\tau _1^p\simeq{c'\tau _1^p},\quad {\left\| {{{\mathbf{b}}_{ex}} - {{\mathbf{b}}_2}} \right\|_b} = c_{2}'\tau _2^p\simeq{c'\tau _2^p},\quad {\left\| {{{\mathbf{b}}_{ex}} - {{\mathbf{b}}_3}} \right\|_b} = c_{3}'\tau _3^p\simeq{c'\tau _3^p},\]
thus
\[{{\mathbf{b}}_{ex}} \simeq \frac{{{{\left\| {\tau _2^p{{\mathbf{b}}_1} - \tau _1^p{{\mathbf{b}}_2}} \right\|}_b}}}{{\tau _2^p - \tau _1^p}},\quad {{\mathbf{b}}_{ex}} \simeq \frac{{{{\left\| {\tau _3^p{{\mathbf{b}}_2} - \tau _2^p{{\mathbf{b}}_3}} \right\|}_b}}}{{\tau _3^p - \tau _2^p}},\]
according to \eqref{propo} we have
\[\frac{{{{\left\| {{{\mathbf{b}}_2} - {{\mathbf{b}}_1}} \right\|}_b}}}{{{{\left\| {{{\mathbf{b}}_3} - {{\mathbf{b}}_2}} \right\|}_b}}} \simeq {\left( {\frac{{{\tau _1}}}{{{\tau _2}}}} \right)^p}.\]
Hence
\[p \simeq \frac{{\log \left( {\frac{{{{\left\| {{{\mathbf{b}}_2} - {{\mathbf{b}}_1}} \right\|}_b}}}{{{{\left\| {{{\mathbf{b}}_3} - {{\mathbf{b}}_2}} \right\|}_b}}}} \right)}}{{\log \left( {\frac{{{\tau _1}}}{{{\tau _2}}}} \right)}}.\qed\]
\end{proof}
In the following section, the numerical errors are computed based on assumptions described in Lemma~\ref{lemma1}~and~\ref{lemma2}.

\section{Numerical results and discussion}\label{sec6}
On the basis of the described procedure, some problems are solved to illustrate the efficiency and the accuracy of the proposed method.
In the first example, a simple two-dimensional fractional heat-like equation is considered for two different conditions. In the second example, a nonlinear two-dimensional fractional wave-like equation is tested. In the third and fourth problems, two linear TFCDEs are solved to test the impact of external force ($g$) and the final time on the convergence rate of the method. In the fifth and sixth test problems, multi-order time-fractional diffusion-wave equations in bounded homogeneous anisotropic plane bodies are solved. The condition number of system \ref{matrixform} is examined for each example. Since, the size of the matrix $\pmb{A}$ depends on the number of internal points, $M$, and the degree of COM, $K$, the condition number of $\pmb{A}$ can be compared versus $K$ and the length of distance between nodal points. However, most domains are not discretized uniformly. In this regard, suppose $r$ denotes the mean length of all the distances between the internal points and their adjacent points (e.g. see $\textstyle{x_{IP}}$ in Figure~\ref{fig: discrete}). Thus, the numerical results show that the condition number behaves as $Cond\pmb{A}\simeq{r^{-2}(K+1)^{2}}$ for example~\ref{exam1}, and  $Cond\pmb{A}\simeq{r^{-2}(K+1)^{3}}$ for other examples. 
\\
\begin{exm}\label{exam1}
Consider the following two-dimensional time fractional heat-like equation: 
\[D_c^{{\alpha}}u = {u_{xx}} + {u_{yy}}, \quad 0 < {\alpha } \leqslant 1,\quad t > 0\]
subjected to different initial conditions with different domains \cite{momani,shirazi}:
\begin{equation*}
  \begin{split}
(I)\quad \quad \begin{gathered}
  u(x,y,0) = \sin x\sin y,   \hfill \\ 
    \quad 0 < x,y < 2\pi ,  \hfill \\
\end{gathered} 
  \end{split}
\quad\quad
  \begin{split}
  (II)\quad \quad \begin{gathered}
  u(x,y,0) = \cos(\frac{\pi}{2}x) \cos(\frac{\pi}{2}y).   \hfill \\ 
    \quad 0 < x,y < 1 ,  \hfill \\
\end{gathered} 
  \end{split}
\end{equation*}
Here, boundary conditions satisfy the exact solutions:
\begin{equation}
\label{Exact1}
  \begin{split}
  (I)\quad u(x,y,t) = {E_\alpha }\left( { - 2{t^\alpha }} \right)\sin x\sin y,
 \hfill \\
  (II)\quad u(x,y,t) = {E_\alpha }\left( { -\frac{1}{2}{\pi}^2t^\alpha } \right)\cos(\frac{\pi}{2}x) \cos(\frac{\pi}{2}y),
  \end{split}  
\end{equation}
where the following one-parameter Mittag-Leffler function is defined as
\begin{equation*}
{E_\alpha }(z) = \sum\limits_{k = 0}^\infty  {\frac{{{z^k}}}{{\Gamma (\alpha k + 1)}}} ,\quad \alpha  > 0.
\end{equation*}
Figure \ref{error} demonstrates ${L_{\infty}}$ errors and $MRE$ versus the degree $K$ (right) for Example~\ref{exam1} (case \textit{I}) when $N=160$, and $\alpha=0.5$ at $t=1.5$. The convergence rate of COM is estimated that $P_{\tau}$-order$>5$. Furthermore, the condition numbers of $\pmb{A}$, on the figure  \ref{error} (left), are shown versus the polynomial degree, $K$, and the mean length of discrete elements, $r$. It can be numerically deduced that the condition number behaves as $Cond(\pmb{A})\simeq{r^{-2}(K+1)^{2}}$; e.g. when $K=8$, then $Cond(\pmb{A})\simeq{76.4\times r^{-2}}$, and when $r=0.3$, then $Cond(\pmb{A})\simeq{10.43\times K^2}$.
In Table \ref{table1}, numerical results are compared with the exact solutions \eqref{Exact1} for Example~\ref{exam1}, case $I$, for fixed $K=12$, $t=1.5$, with the differential orders $\alpha=0.5$ and $\alpha=0.75$, and different number of nodal points; the convergence rate of BEM is algebraic ($P_{r}$-order$>4.4$) when the number of nodal points is increased from $N=40$ to $N=80$ and it is quadratic ($P_{r}$-order$>2$) when $N=80$ goes to $N=160$. Apart from the value of $\alpha$, it can be inferred that the computation cost of the second discretization for moderate $N$ and $M$ is more effective than the third one.\\
For case ({\it{II}}), a similar behavior of $Cond(\pmb{A})$ versus $K$ and $r$ is shown in Figure \ref{figProb1case2} (left). The relative absolute error (right) with $N=200$, and $K=12$ for the final time $t=1$ is exhibited. Intuitively, the relative absolute errors are approached to $10^{-5}$, which it could be expected for $K=12$ and $N=200$ based on the information from~\ref{table1II}. This table  shows the estimated convergence for two terms of shifted Chebyshev series for the final time $t=0.5$; in general, there is an improvement for errors when the degree $K$ increases, but no relationship between $P_r$-order and degree $K$ is observed. It may also be concluded $N=64$ is more computational cost-effective in this case.

\begin{figure}[!ht]
\begin{center}
\subfigure{\includegraphics[width=.6\textwidth]{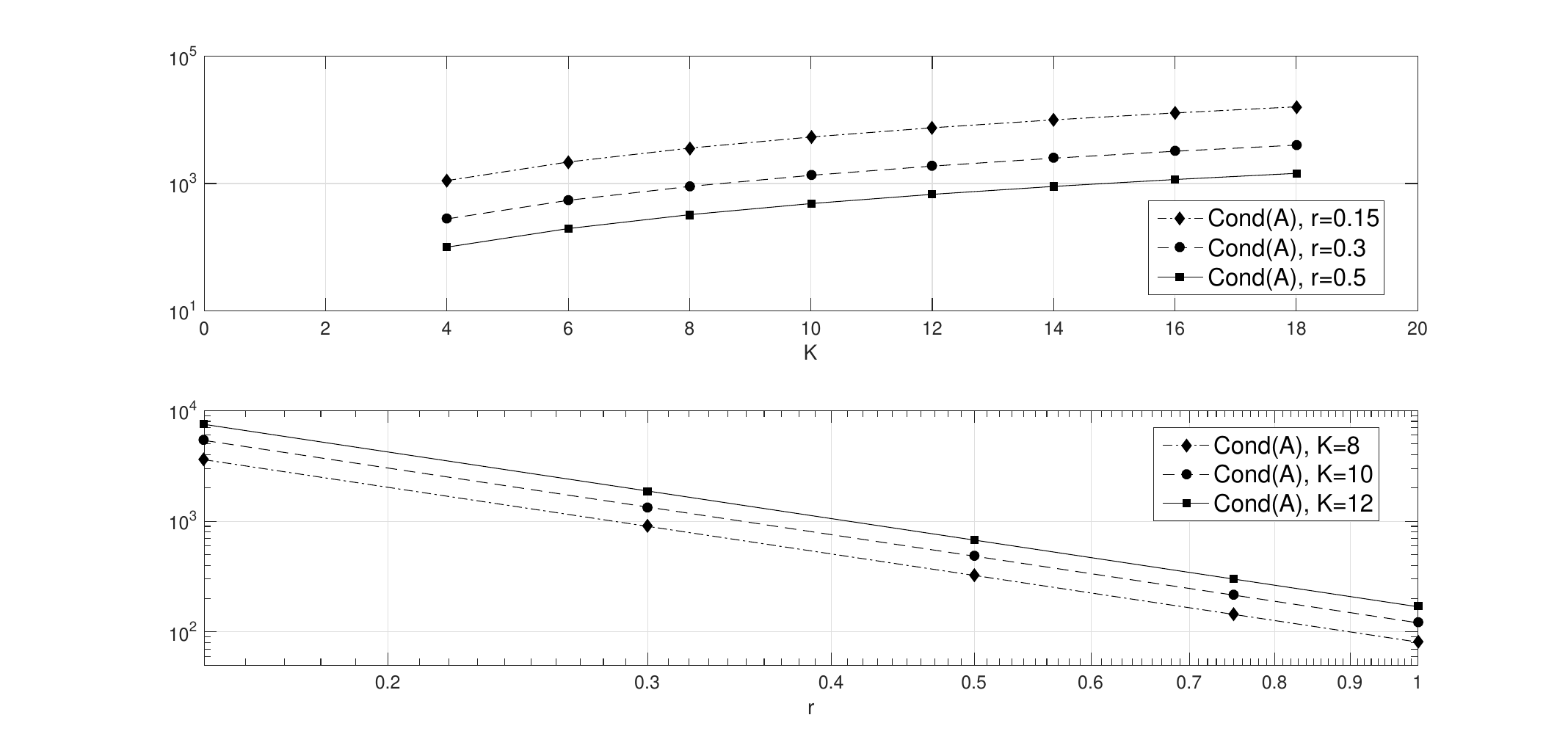}}
\subfigure{\includegraphics[width=.38\textwidth]{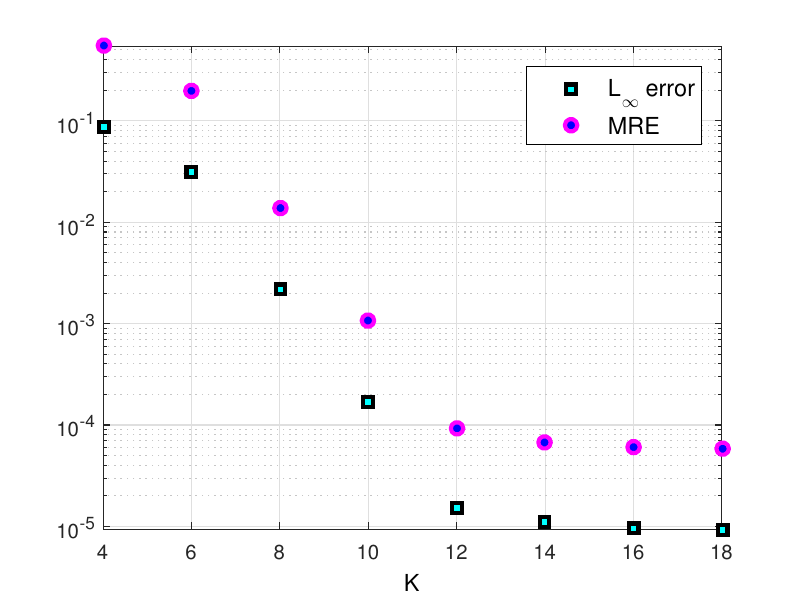}}
\caption{The condition number of $\pmb{A}$ versus the polynomial degree $K$ and the mean length of discrete elements $r$ (left) and a comparison of $L_\infty$ errors and $MRE$ versus $K$ (right) at  $t=1.5$ when $N=160$ and $\alpha=0.5$ for Example~\ref{exam1} case $I$; $P_{\tau}$-order$=5.13$ is estimated for $K=4, 8, 16$.}\label{error}
\end{center}
\end{figure}

\begin{table}[!ht]
\caption{The error norms and the estimated order of convergence $P_{r}$ for the vector solution ${\mathbf{U}}$ according to the Lemma~\ref{lemma1}, in Example~\ref{exam1} case $I$ when $K=12$ and $t=1.5$.} 
\centering 
\begin{tabular}{|c| c c c| c c c|} 
\hline
\multirow{2}{*}{$N$} &\multicolumn{3}{c|}{$\alpha=0.5$} &\multicolumn{3}{c|}{$\alpha=0.75$}\\
\cline{2-7}
  & $MRE$ & $RMS$ & $P_{r}$-order & $MRE$ & $RMS$ & $P_{r}$-order\\ [1ex]
\hline 
 40 & $5.48066\times {10^{-3}}$ & $2.38289\times {10^{-4}}$ & $-$ & $8.21377\times {10^{-3}}$ & $3.57120\times {10^{-4}}$  & $-$\\ 
80 & $2.52922\times {10^{-4}}$ & $1.09966\times {10^{-5}}$ & 4.4376 & $3.72158\times {10^{-4}}$ & $1.61807\times {10^{-5}}$ & 4.4641\\ 
160 & $6.10497\times {10^{-5}} $ & $2.65433\times {10^{-6}}$ & 2.0506 &$8.43489\times {10^{-5}}$ & $3.66734\times {10^{-6}}$ & 2.1415\\ 
\hline
\end{tabular}
\label{table1}
\end{table}

\begin{figure}[!ht]
\begin{center}
\subfigure{\includegraphics[width=.6\textwidth]{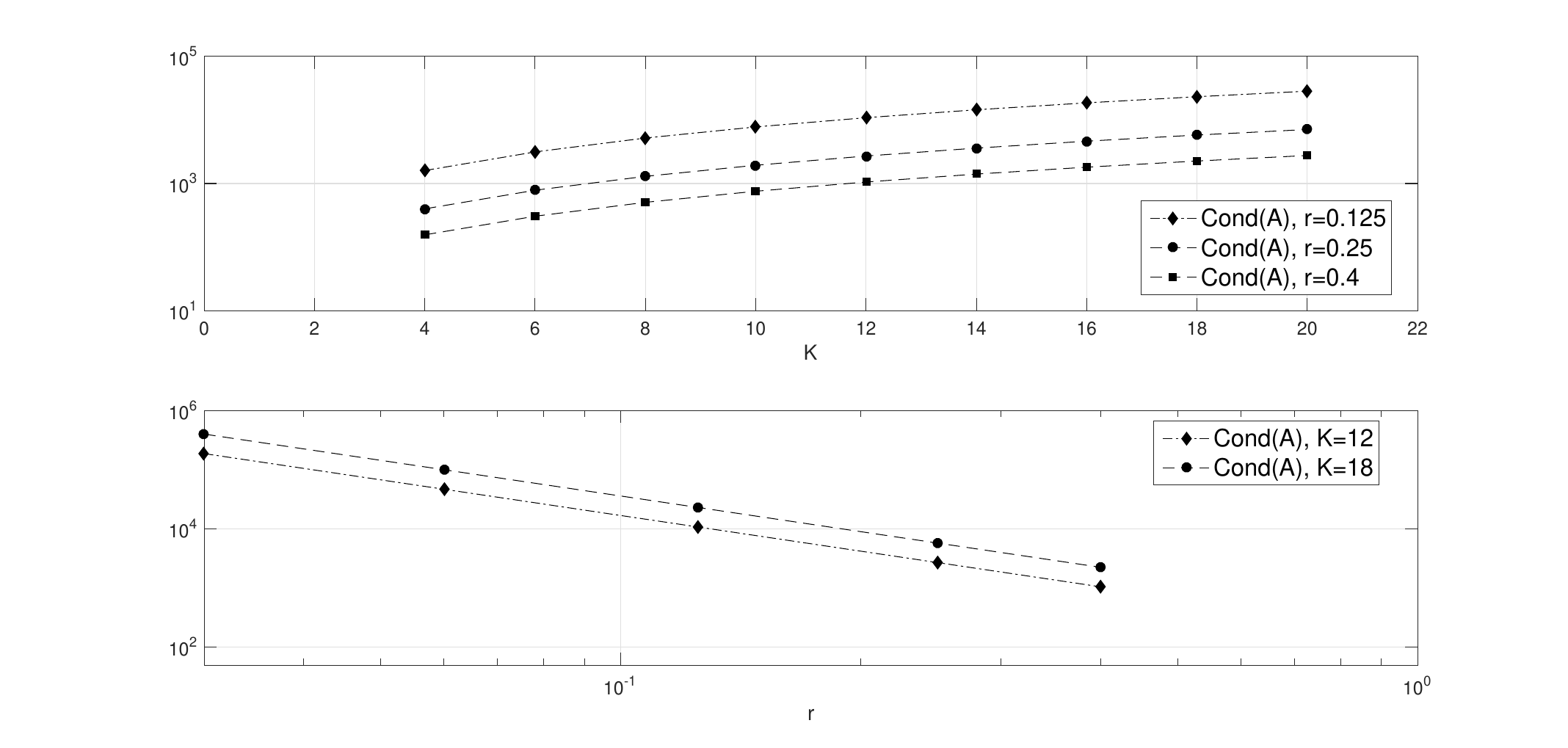}}
\subfigure{\includegraphics[width=.38\textwidth]{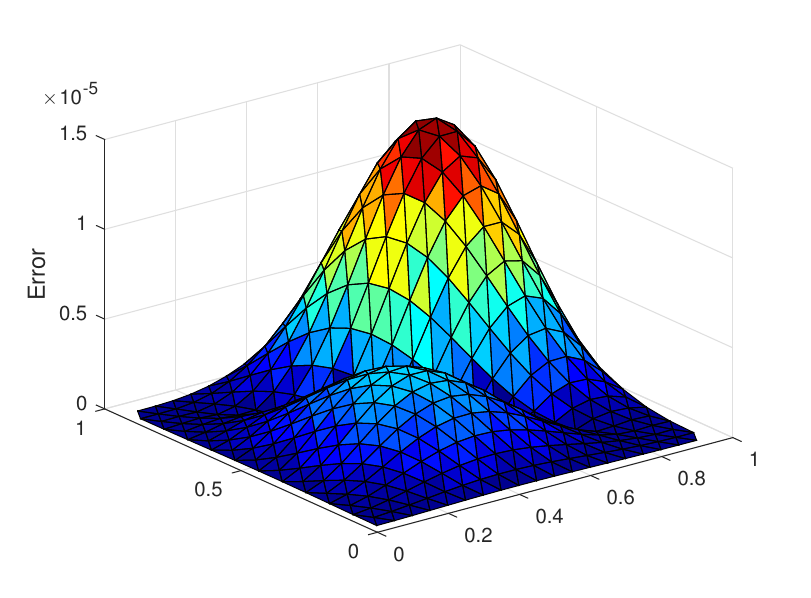}}
\caption{The condition number of $\pmb{A}$ versus the polynomial degree $K$ and the mean length of discrete elements $r$ (left) and the relative absolute errors (right) obtained for Example~\ref{exam1} case $II$ when $N=200$, $K=12$, $\alpha=0.5$, and $t=0.5$.}\label{figProb1case2}
\end{center}
\end{figure}

\begin{table}[!ht]
\begin{center}
\caption{A comparison of the error norms and the estimated order of convergence $P_{r}$ for the vector solution ${\mathbf{U}}$ according to the Lemma~\ref{lemma1}, in Example~\ref{exam1} case $II$ for two fixed $K$.} 
\begin{tabular}{|c| c c c| c c c|} 
\hline
 \multirow{2}{*}{$N$} & \multicolumn{3}{c|}{$K=12$} &\multicolumn{3}{c|}{$K=18$}\\
 \cline{2-7}
  & $L_\infty$ & $RMS$ & $P_{r}$-order & $L_\infty$ & $RMS$ & $P_{r}$-order\\ [1ex]
\hline 
 8 & $1.9970\times {10^{-3}}$ & $1.8930\times {10^{-3}}$  & $-$ & $1.3222\times {10^{-3}}$ & $1.0553\times {10^{-3}}$ & $-$\\ 
16 & $1.4673\times {10^{-3}}$ & $7.1121\times {10^{-4}}$ & $1.4123$ & $1.0571\times {10^{-3}}$ & $5.1238\times {10^{-4}}$  & $1.0424$\\ 
 32 & $6.0501\times {10^{-4}} $ & $1.2220\times {10^{-4}}$ & $2.5410$ &$3.5453\times {10^{-4}}$ & $7.1607\times {10^{-5}}$  & $2.8390$\\ 
 64 & $9.3468\times {10^{-5}}$ & $9.9095\times {10^{-6}}$ & $3.6243$ & $6.9059\times {10^{-5}}$  & $7.3216\times {10^{-6}}$  & $3.2899$\\ 
 128 & $ 3.8455\times {10^{-5}} $  & $1.8098\times {10^{-6}}$ & $2.4530$ & $2.4404\times {10^{-5}}$ & $1.1486\times {10^{-6}}$ & $2.6723$\\               	\hline
\end{tabular}
\label{table1II}
\end{center}
\end{table}
\end{exm}

\begin{exm}\label{exam2}
Consider the two-dimensional time fractional wave-like equation \cite{sho}: 
\[D_c^{{\alpha }}u = \frac{1}{12}\left( {{x^2}{u_{xx}} + {y^2}{u_{yy}}} \right),\quad 0 < x,y < 1,\quad 1 < {\alpha} \leqslant 2,\quad t > 0,\]
subjected to boundary conditions 
\[\begin{gathered}
  u(0,y,t) = 0,\quad u(1,y,t) = 4\cosh t, \hfill \\
  u(x,0,t) = 0,\quad u(x,1,t) = 4\sinh t ,\hfill \\ 
\end{gathered} \]
and the initial condition 
\[\begin{gathered}
  u({\mathbf{x}},0) = {x^4}, \hfill \\
  {u'}({\mathbf{x}},0) = {y^4}. \hfill \\
\end{gathered} \]
The exact solution for $\alpha=2$ is found to be,
\begin{equation}
\label{Exact2}
u(x,y,t) = {x^4}\cosh t + {y^4}\sinh t.
\end{equation}
This problem is solved for different $K$ with $N=100$ at $t=0.5$ and integer order $\alpha=2$ to compared with the exact solution \eqref{Exact2}. The results from Table \ref{table2} offer an error improvement by increasing the number of degree $K$ of COMs. The condition numbers of the system \ref{matrixform}  illustrate such behaviour $Cond(\pmb{A})\simeq{r^{-2}(K+1)^{3}}$ (see Table \ref{table2}). 
Due to the fact that the domain and boundary nodal points are fixed here, the numerical solutions of $\textbf{U}$ are directly affected by the numerical solutions of $\textbf{b}$; in other words, affected by the accuracy of COM. However, the exact solution of the generated ODE system is not clear, and the convergence order of COMs is estimated by norm ${\left\| . \right\|_b}$ for three distinct degrees with a same proportion. Interestingly, a comparison between the two columns $RMS$ and ${\left\| . \right\|_b}$ of  Table~\ref{table2} suggests direct relationships, but with different speed between the approximation solution of $\textbf{U}$ and $\textbf{b}$. In addition, by considering the scale of the solutions, and a comprehensive assessment of the absolute errors for distinct degrees $K$, it could be concluded that the Chebyshev Tau method converges with an oscillating manner around the exact solution of fractional ODE system~\eqref{eq14}.

\begin{table}[!ht]
\caption{The error norms for the vector solution ${\mathbf{U}}$  and the convergence orders of the vector solution ${\mathbf{b}}$ in Example~\ref{exam2}, for $\alpha=2$, $N=100$ and $t=0.5$.} 
\vspace{.3cm}	
\centering 
\begin{tabular}{|c c c c c c c|} 
\hline
 $K$ & $L_\infty$ &  $MRE$ & $RMS$ & ${\left\| {{{\mathbf{b}}_i} - {{\mathbf{b}}_{i - 1}}} \right\|_b}$ & $P_{\tau}$-orde & $Cond(\pmb{A})$\\ [1ex]
\hline 
 8  & $7.15174\times {10^{-3}}$ & $2.25751\times {10^{-2}}$ & $2.91855\times {10^{-3}}$  & $-$ & $-$ & $4.54024 \times {10^{5}}$\\ 
16  & $2.82029\times {10^{-5}}$ & $8.43598 \times {10^{-5}}$ & $1.06609\times {10^{-5}}$ & $5.5885 \times {10^{ - 6}}$& $-$ & $2.96189 \times {10^{6}}$\\ 
32 & $1.07633\times {10^{-5}}$ & $2.08835\times {10^{-5}}$ & $6.63201\times {10^{-6}}$ & $5.0329\times {10^{-7}}$  & $3.4730$ & $2.14632 \times {10^{7}}$\\
64 & $8.17112\times {10^{-6}}$ & $5.69019\times {10^{-6}}$ & $9.87053\times {10^{-7}}$ & $3.2380\times {10^{-8}}$  & $3.9582$ & $1.71381 \times {10^{8}}$\\
\hline
\end{tabular}
\label{table2}
\end{table}

\end{exm}

\begin{exm}\label{exam3}
Consider the following TFCDE \cite{wang2014}:
\[\begin{gathered}
  D_c^\alpha u = {u_{xx}} + {u_{yy}} - 5{u_x} - 5{u_y} + g(x,y,t), \hfill \\
  \quad \quad 1 < \alpha  \leqslant 2,\quad (x,y) \in \Omega , \hfill \\ 
\end{gathered} \]
with the boundary condition and initial conditions
\[\begin{gathered}
  u(x,y,t) = 0,\quad (x,y) \in \Gamma ,\quad 0 < t \leqslant 1, \hfill \\
  u(x,y,0) = 0,\quad \psi (x,y) = 0,\quad (x,y) \in \Omega , \hfill \\ 
\end{gathered} \]
where $\Omega  = {[0,1]^2}$. Then we have the following exact solution:
\[u(x,y,t) = {2^{12}}{t^{2 + \alpha }}{x^3}{(1 - x)^3}{y^3}{(1 - y)^3}.\]
It is easy to check
\[\begin{gathered}
  g(x,y,t) = {2^{11}}\Gamma (\alpha  + 3){x^3}{(1 - x)^3}{y^3}{(1 - y)^3}{t^2} \hfill \\
  \quad  - {2^{12}}(6x - 51{x^2} + 120{x^3} - 150{x^3} - 150{x^4} + 30{x^5}){y^3}{(1 - y)^3}{t^{2 + \alpha }} \hfill \\
  \quad  - {2^{12}}{x^3}{(1 - x)^3}(6y - 51{y^2} + 120{y^3} - 105{y^4} + 30{y^5}){t^{2 + \alpha }}. \hfill \\ 
\end{gathered} \]
This problem is challenging, and sensitive because of the large numbers included in the function $g$. However, it could be compensated by multiplying to power functions of decimal numbers, and considering the final time $t=1$. Figure \ref{figProb3} (left plan) shows the estimated error ranged around $10^{-4}$, for $\alpha=1.5$ and final time $t=1$, with the degree $K=10$, and (right plan) demonstrates the plot of the error versus the number of boundary nodes, $N$, with $K=10$ for three different values $\alpha$, illustrating that the smoothness roughly occurred after $N=135$.
The behavior of the condition number matrix $\pmb{A}$ is estimated as $Cond(\pmb{A})\simeq{r^{-2}(K+1)^{3}}$. 
Table~\ref{tableEx321} gives $Cond(\pmb{A})$, the RMS error and the convergence rates are obtained by solving Example~\ref{exam3} for different values of $\alpha$. It indicates better $RMS$ errors for $\alpha$ near to 1 than 2, that is not true for $P_r$-orders.

\begin{figure}[!ht]
\begin{center}
\subfigure{\includegraphics[width=.4\textwidth]{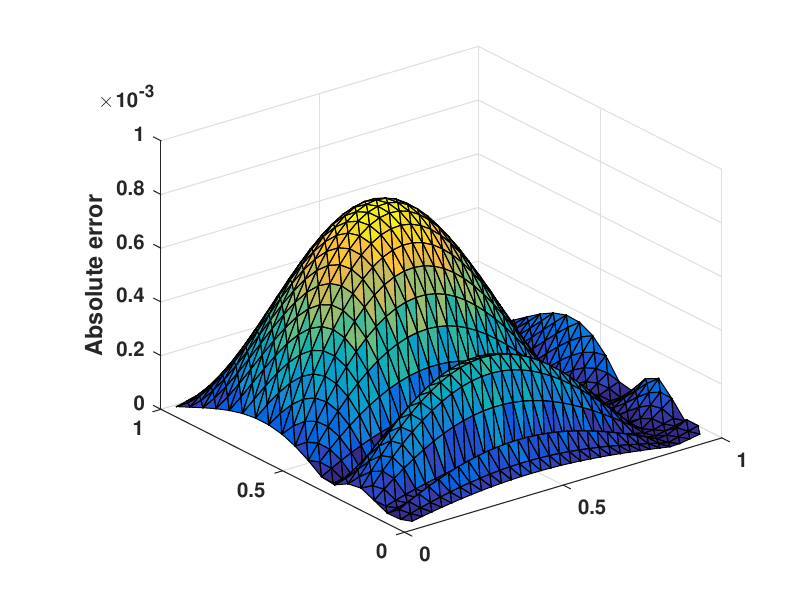}}
\subfigure{\includegraphics[width=.4\textwidth]{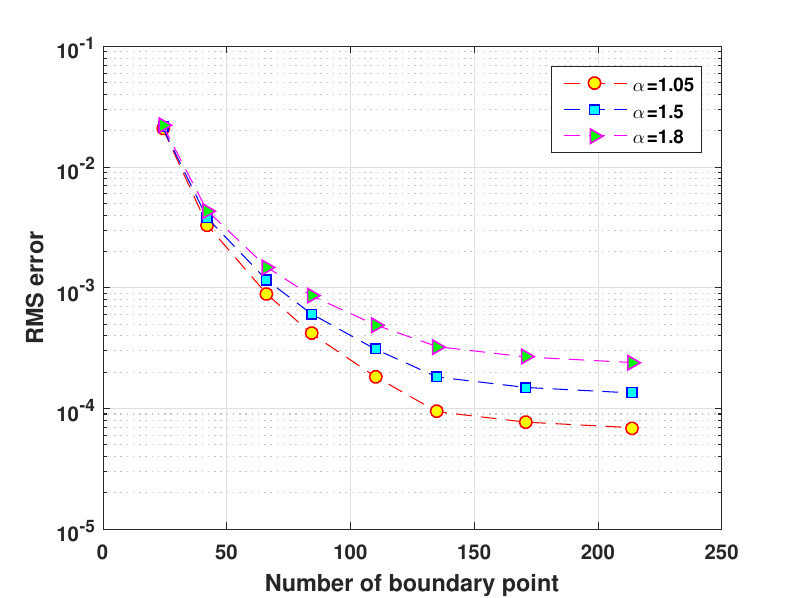}}
\caption{The graph of absolute errors with $\alpha=1.5$ and $K=10$ (left) and a comparison of errors for different values of $\alpha$ (right) at finite time $t=1$ for Example~\ref{exam3}.}\label{figProb3}
\end{center}
\end{figure}

\begin{table}[!ht]
\caption{The $RMS$ error of the vector solution ${\mathbf{U}}$ and the convergence rate of COM and spatial-discretization base on Lemmas~\ref{lemma1} and ~\ref{lemma2} , for Example~\ref{exam3} at  time $t=1$.} 
\vspace{.3cm}
\centering 

\begin{tabular}{|c| c c c | c c c|} 
\hline 
$K=16$ & \multicolumn{3}{c}{$\alpha=1.9$} & \multicolumn{3}{c|}{$\alpha=1.1$} \\
\hline 
$N$ & $RMS$ & $Cond(\pmb{A})$ & $P_{r}$-order & $RMS$ & $Cond(\pmb{A})$ &  $P_{r}$-order\\ [1ex]
\hline 
  40 &$8.8704\times {10^{-3}}$  & $4.86686 \times {10^{5}} $ & $-$ & $7.3330\times {10^{-3}}$ & $ 4.89036 \times {10^{5}}$ & $-$ \\ 
 80   &$2.0840\times {10^{-3}}$  & $ 1.95421 \times {10^{6}}$ & 2.0897& $1.6591\times {10^{-3}}$ & $ 1.96064 \times {10^{6}}$ & 2.1439\\ 
 160 &$4.8292\times {10^{-4}}$ & $7.78881 \times {10^{6}}$ & 2.1095& $3.5514\times {10^{-4}}$  & $ 7.80051 \times {10^{6}}$ & 2.2240 \\ 
  320  & $8.1966\times {10^{-5}}$ & $ 3.13058 \times {10^{7}}$ & 2.5587& $    7.2452\times {10^{-5}} $  & $ 3.13134 \times {10^{7}}$ & 2.2933\\       

\noalign{\hrule height 2pt}
 $N=200$ & \multicolumn{3}{c}{$\alpha=1.9$ }&\multicolumn{3}{c|}{$\alpha=1.6$ }\\
\hline
 $K$ & $RMS$ & ${\left\| {{{\mathbf{b}}_i} - {{\mathbf{b}}_{i - 1}}} \right\|_b}$  & $P_{\tau}$-order & $RMS$ & ${\left\| {{{\mathbf{b}}_i} - {{\mathbf{b}}_{i - 1}}} \right\|_b}$  & $P_{\tau}$-order\\ [1ex]
\hline 
 10 & $1.3010 \times {10^{-4}}$ & $-$ & $-$ & $7.5314\times {10^{-5}}$  & $-$ & $-$ \\ 
20 &  $  4.2187\times {10^{-5}}$ & $1.3840 \times {10^{-6}}$ & $-$ & $1.6938\times {10^{-5}}$ & $1.4839\times {10^{-6}}$ & $ -$\\ 
40 &  $1.2781\times {10^{-5}}$ & $  3.7002\times {10^{-7}}$ &  $  1.9032$ & $3.8515\times {10^{-6}}$  & $2.8012\times {10^{-7}}$  & $ 2.4053$ \\ 
 80 & $ 3.7745\times {10^{-6}}$ &  $ 9.8207\times {10^{-8}}$  & $  1.9137$ & $8.8808\times {10^{-7}}$  &  $4.8922\times {10^{-8}}$ & $   2.5175$ \\       
    \hline
\end{tabular}
\label{tableEx321}
\end{table}
\end{exm}

\begin{exm}\label{exam4}
Consider the linear TFCDE \cite{dehghan2016}:
\[\begin{gathered}
  D_c^\alpha u = {u_{xx}} + {u_{yy}} - 0.1{u_x} - 0.1{u_y} + g(x,y,t), \hfill \\
  \quad \quad 1 < \alpha  \leqslant 2,\quad (x,y) \in \Omega , \hfill \\ 
\end{gathered} \]
with the initial conditions
\[ u(x,y,0)=0,\; u'(x,y,0)=0, \quad (x,y) \in \Omega , \]
and the Dirichlet boundary conditions. The exact solution of the current test problem is
\[u(x,y,t) = {t^{3 + \alpha }}\sin \left( {\frac{\pi }{6}x} \right)\sin \left( {\frac{{7\pi }}{4}x} \right)\sin \left( {\frac{{3\pi }}{4}y} \right)\sin \left( {\frac{{5\pi }}{4}y} \right),\]
where $\Omega$ is the computational domain as shown in figure \ref{fig1Ex4} (left plan). The approximate solution and its relative absolute error are shown in figure \ref{fig1Ex4} for $\alpha=1.05$ and $\alpha=1.95$ with $K=10$ and $N=255$ for the final time $t=2$. Although $N=255$ is considered to depict clear data points in the figure, $N=150$ would be sufficient to achieve the semi-equivalent errors. Importantly, by considering the results of the previous example~\ref{exam3} and Figure~\ref{fig1Ex4} (right plan), it may convey that the method has a better performance for the less values of $\alpha$. In contrast, Table~\ref{table4} refuses this idea; there are irrelevant outcomes versus the values of $\alpha$, although the table shows a reliable numerical convergence.

\begin{figure}[!ht]
\begin{center}
\subfigure{\includegraphics[width=.4\textwidth]{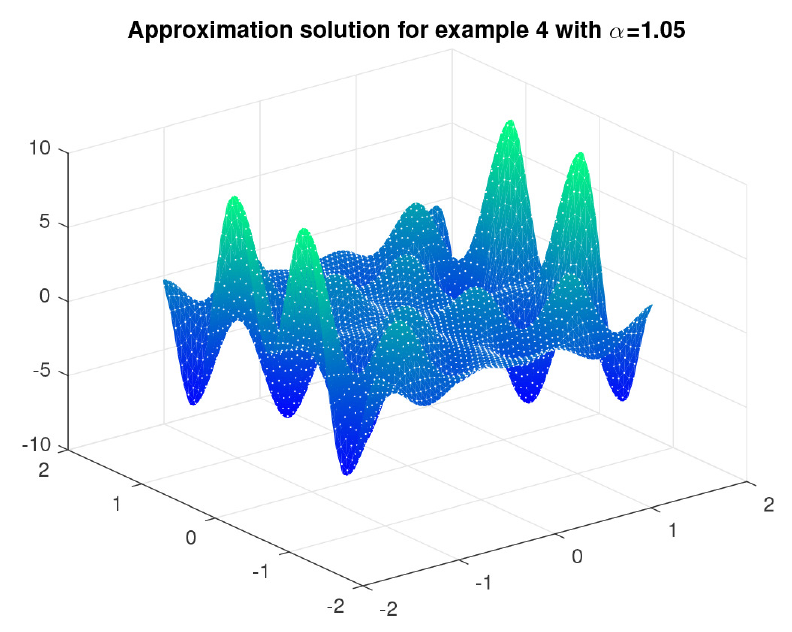}}
\subfigure{\includegraphics[width=.4\textwidth]{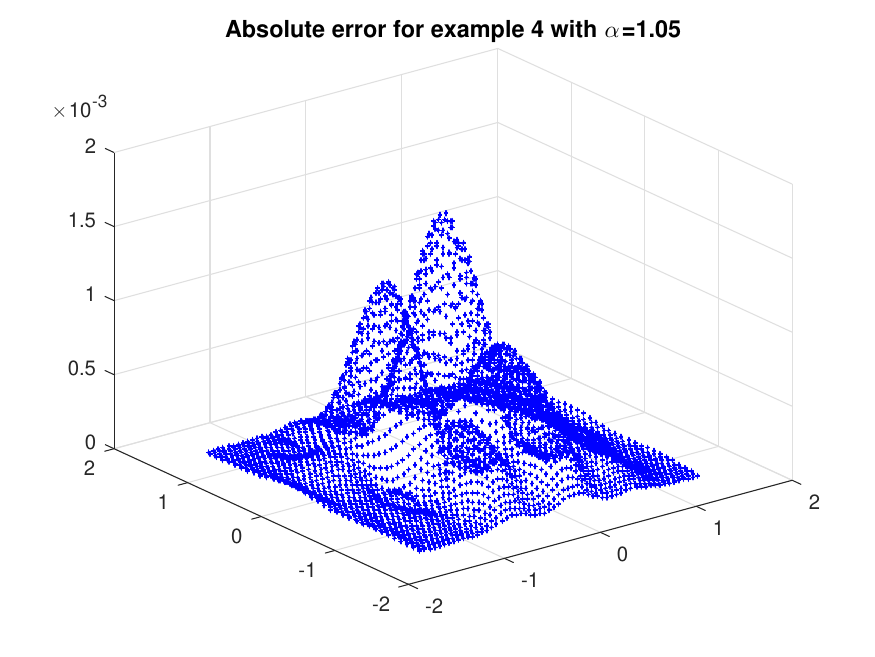}}
\subfigure{\includegraphics[width=.4\textwidth]{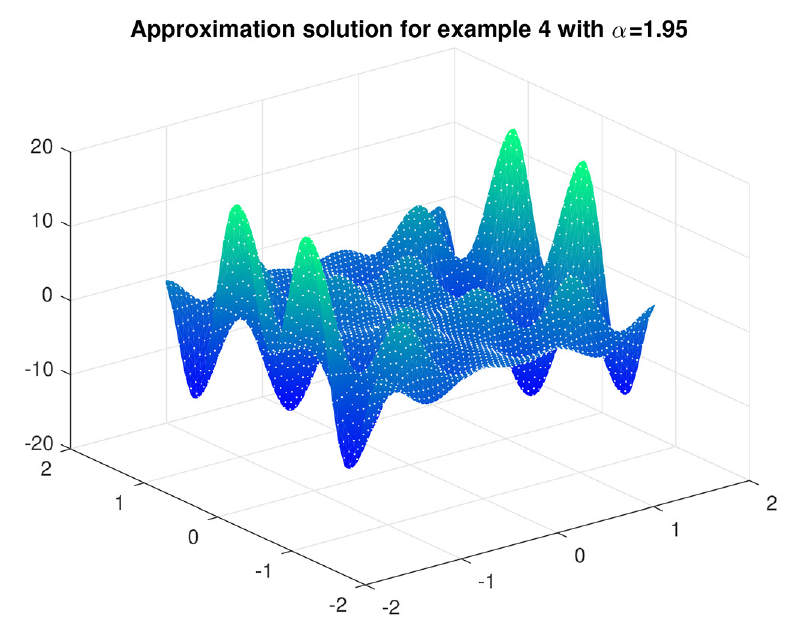}}
\subfigure{\includegraphics[width=.4\textwidth]{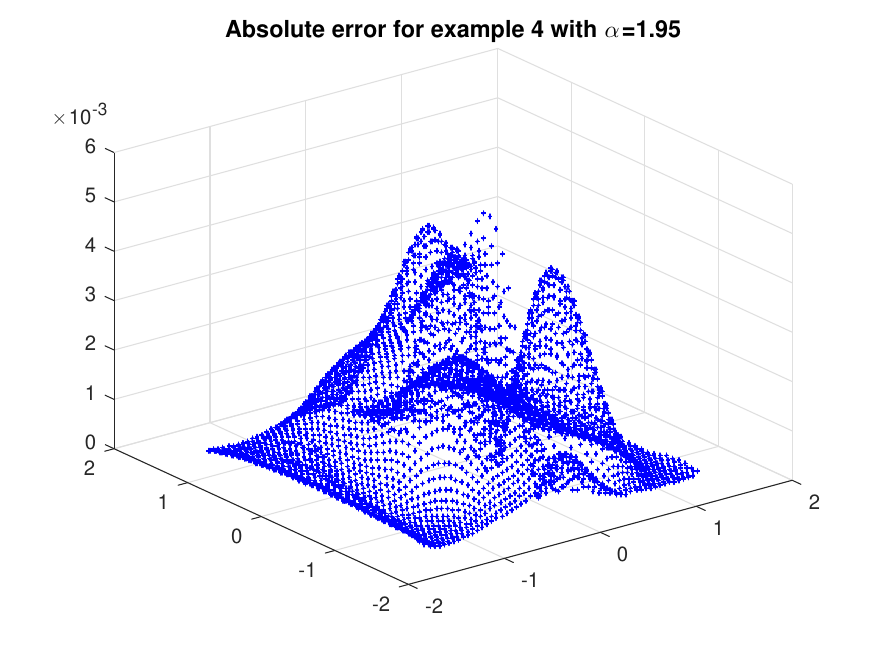}}
\caption{Graphs of approximation solution and the absolute error obtained for Example~\ref{exam4} for $\alpha=1.05$ and $\alpha=1.95$ with $K=10$ and $N=255$ for final time $t=2$.}\label{fig1Ex4}
\end{center}
\end{figure}

\begin{table}[!ht]
\caption{The $RMS$ error for the vector solution ${\mathbf{U}}$, in Example~\ref{exam4} for $K=8$ and $t=2$.} 
\vspace{.3cm}
\centering 
\begin{tabular}{|c c c c c|} 
\hline 
\multicolumn{5}{|c|}{$RMS-error$} \\
\hline 
$N$  & $\alpha=1.1$ & $\alpha=1.4$  &  $\alpha=1.7$ & $\alpha=1.9$ \\ [1ex]
\hline 
  10 & $1.5030\times {10^{-1}}$ & $1.0628\times {10^{-1}}$ & $1.2802\times {10^{-1}}$ & $1.1669\times {10^{-1}}$ \\ 
 20  & $3.5223\times {10^{-2}} $ & $3.8318\times {10^{-2}}$ & $5.0480\times {10^{-2}}$ & $4.7404\times {10^{-2}}$\\ 
 40 & $   1.6236\times {10^{-2}}$ & $1.1808\times {10^{-2}}$ & $1.4409\times {10^{-2}}$  & $1.1033\times {10^{-2}}$ \\ 
 80 & $    2.4654\times {10^{-3}} $  & $4.9065\times {10^{-3}}$ &  $2.0146\times {10^{-3}}$ & $4.5839\times {10^{-3}}$ \\             160 & $   1.2483\times {10^{-3}} $  & $9.5059\times {10^{-4}}$ &  $9.5252\times {10^{-4}}$ & $1.1706\times {10^{-3}}$ \\                     320 & $   3.7848\times {10^{-4}} $  & $2.0804\times {10^{-4}}$ &  $1.7544\times {10^{-4}}$ & $2.1400\times {10^{-4}}$ \\               \hline 
\end{tabular}
\label{table4}
\end{table}
\end{exm}

\begin{exm}\label{exam5}
The multi-order time fractional diffusion-wave equation \cite{katsi2011}
\begin{equation}
\label{diffusion}
{\gamma _1}D_c^{1.7}u + {\gamma _0}D_c^{0.8}u = A{u_{xx}} + 2B{u_{xy}} + C{u_{yy}} + g({\mathbf{x}},t),\quad {\mathbf{x}}(x,y) \in \Omega ,\quad t > 0,
\end{equation}
in the plane inhomogeneous anisotropic body which is shown in Figure. \ref{fig1} has been solved, subject to boundary conditions 
\[u({\mathbf{x}},t) = 0,\quad {\mathbf{x}}(x,y) \in \Gamma ,\]
and the initial condition 
\[\begin{gathered}
  u({\mathbf{x}},0) = 0, \hfill \\
  {u'}({\mathbf{x}},0) = U(x,y), \hfill \\ 
\end{gathered} \]
where
$A = \frac{{({y^2} - {x^2} + 50)}}{{50}}$,
$B = \frac{{2xy}}{{50}}$,
$C = \frac{{({x^2} - {y^2} + 50)}}{{50}}$,
${\gamma _1} = 5{e^{\left( { - 0.1\left( {\left| x \right| + \left| y \right|} \right)} \right)}}$,
${\gamma _0} = 0.4{\left( {{x^2} + {y^2}} \right)^{1/2}}$. The external source $g$ is
$g({\mathbf{x}},t) =U({\gamma _1}D_c^{1.7}T + {\gamma _0}D_c^{0.8}T) -T(A{U_{xx}} + 2B{{U}_{xy}} + C{U_{yy}}),$
where $U(x,y) = {a^2}{b^2} - \left( {{{\left( {\frac{x}{a}} \right)}^2} + {{\left( {\frac{y}{b}} \right)}^2}} \right)\left( {{{\left( {\frac{x}{b}} \right)}^2} + {{\left( {\frac{y}{a}} \right)}^2}} \right)$ and $T(t) = t - \frac{{{t^3}}}{6} + \frac{{{t^5}}}{{200}}$. 
The boundary of the domain is defined by the curve:
\[\Gamma  = \frac{{{{\left( {ab} \right)}^{1/2}}}}{{{{\left( {{{\left( {\cos \frac{\theta }{a}} \right)}^2} + {{\left( {\sin \frac{\theta }{b}} \right)}^2}} \right)}^{(1/4)}}{{\left( {{{\left( {\cos \frac{\theta }{b}} \right)}^2} + {{\left( {\sin \frac{\theta }{a}} \right)}^2}} \right)}^{(1/4)}}}},\]
where $0 \leqslant \theta  \leqslant 2\pi ,\quad a = 3,\quad b = 1.3.$ The problem admits an exact solution ${u_{exact}} =T(t){U}(x,y)$. This problem is solved using BEM and COM for various $M$ and $N$ when $K=16$ at $t=4$. The value of $u$, $u_x$ and $u_{xy}$ are compared with the exact solution. Numerical results are given in Table \ref{table5} showing the efficiency of the proposed method by $P_{r}$-order $>2$, and the condition number of matrix $\pmb{A}$ with the manner as $r^{-2}\times (K+1)^3$. In figure \ref{fig1}, the contour plots illustrate the absolute errors distributions of the approximations of $u$, $u_x$ and $u_{xy}$ on the plane, including $L_{\infty}, MRE, RMS$, at $t=4$ for specific nodal points and Delaunay triangulation when $N=210$, $M=132$, $K=16$. 

\begin{figure}[!ht]
\begin{center}
\subfigure{\includegraphics[width=.35\textwidth]{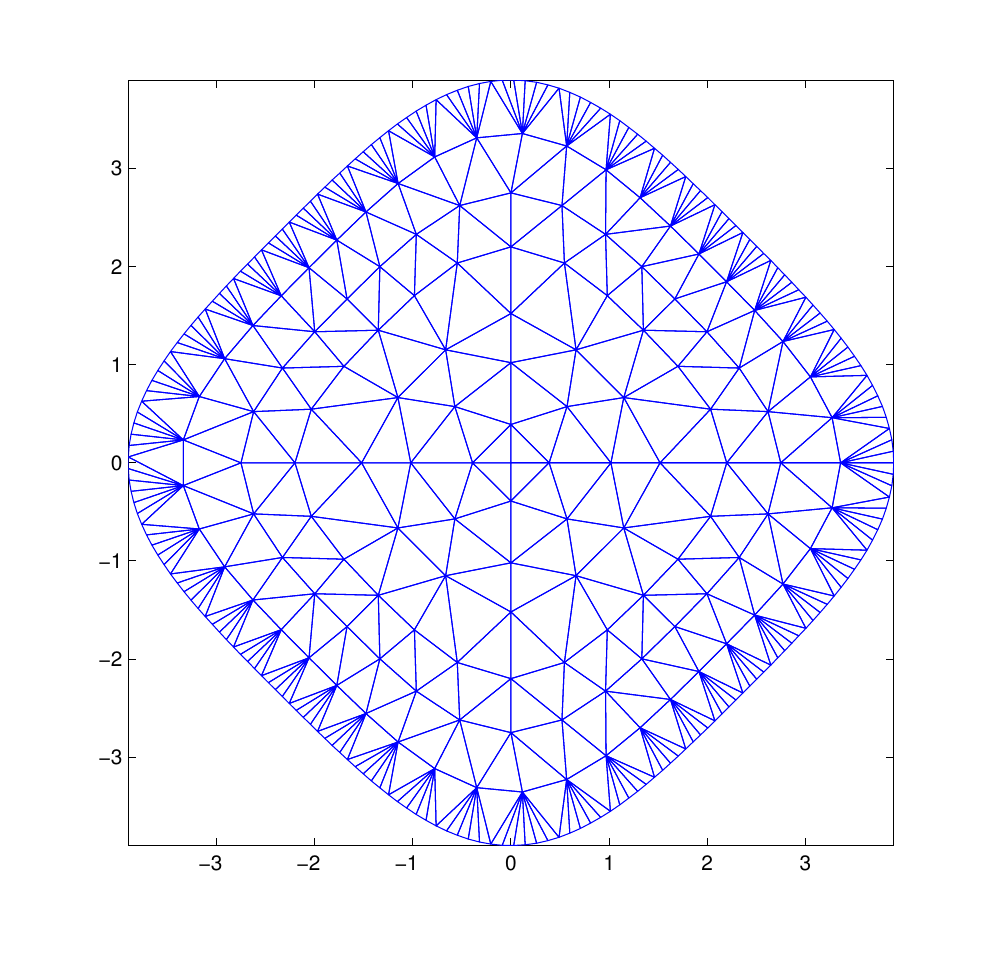}}
\subfigure{\includegraphics[width=.4\textwidth]{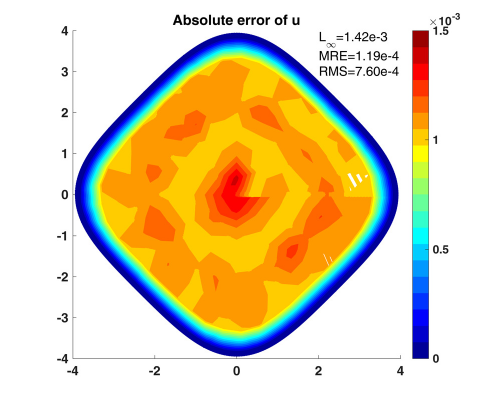}}
\subfigure{\includegraphics[width=.4\textwidth]{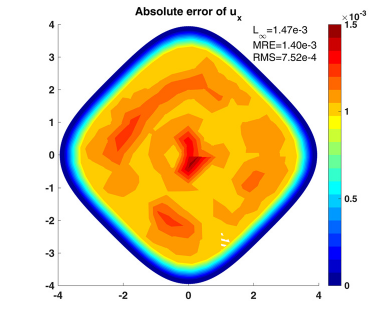}}
\subfigure{\includegraphics[width=.4\textwidth]{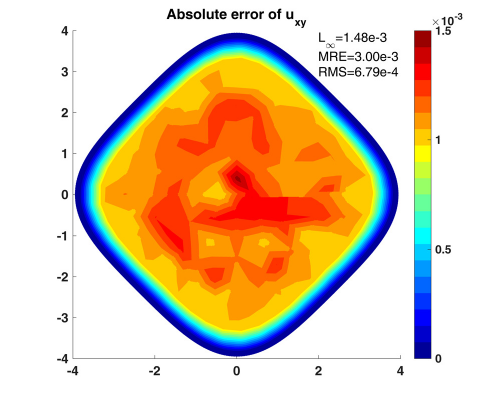}}
\caption{The geometry of the plane body discretization with  $N=210$, $M=132$ with DG-RBF technique and Contour plots of absolute errors when $K=16$ for Example~\ref{exam5}.}\label{fig1}
\end{center}
\end{figure}

\begin{table}[!ht]
\caption{The condition number of $\pmb{A}$, the $RMS$ error and the convergent order for the vector solution ${\mathbf{b}}$ according to the Lemma~\ref{lemma2}, in Example~\ref{exam5} for $K=16$ and $t=4$.} 
\vspace{.3cm}
\centering 
\begin{tabular}{|c c c c c c c|} 
\hline 
\multicolumn{5}{|c}{$RMS-error$}& & \\
\cline{2-4}
$N$  & $u$ & $u_{x}$  &  $u_{xy}$ & $P_{r}$-order & $r\simeq$  & $Cond(\pmb{A})$ \\ [1ex]
\hline
  50 & $6.27668\times {10^{-2}}$ & $6.69339\times {10^{-2}}$ & $6.37693\times {10^{-2}}$ & $-$ & $2.4$ & $8.5723 \times {10^{2}}$ \\ 
 100  & $6.74090\times {10^{-3}} $ & $5.97502\times {10^{-3}}$ & $6.03444\times {10^{-3}}$ & $3.40157$ & $1.2$ & $3.4289 \times {10^{3}}$\\ 
 120 & $  3.70614\times {10^{-3}}$ & $3.60496\times {10^{-3}}$ & $3.86762\times {10^{-3}}$ & $2.43988$ & $0.1$ & $4.9376 \times {10^{3}}$ \\ 
 200 & $  7.71852\times {10^{-4}} $  & $7.85751\times {10^{-4}}$ &  $7.03289\times {10^{-4}}$ & $3.33700$ & $0.6$ & $1.3715 \times {10^{4}}$\\
  \hline 
\end{tabular}
\label{table5}
\end{table}

\end{exm}

\begin{exm}\label{exam6}
Consider the large terms time fractional diffusion equation

\begin{equation}
\label{e2exam6}
\begin{gathered}
  \sum\limits_{j = 0}^6 {{\gamma _j}({\mathbf{x}})D_c^{{\alpha _j}}u}  = D({\mathbf{x}}){u_x} + E({\mathbf{x}}){u_y} + F({\mathbf{x}})u + g({\mathbf{x}},t), \quad {\mathbf{x}}(x,y) \in \Omega ,\quad t > 0, \hfill 
\end{gathered} 
\end{equation}
in a ``\textbf{C}-shape'' made by the elimination of a circle with radius $r_{2}=3$ and null origin from the inside of a circle with radius $r_{1}=5$ and the same origin, and extracting the space between the lines $y=-1$ and $y=1$ from the right side of the outcome (see Figure~\ref{figEX6}) with the Dirichlet boundary conditions, and the initial condition
\[\begin{gathered}
  u({\mathbf{x}},0) = 0, \hfill \\
  {u'}({\mathbf{x}},0) = U(\textbf{x}), \hfill \\ 
\end{gathered} \]
where $U(\textbf{x})=cos(x)e^{sin(y)}$, and the external force $g$ is $g(\textbf{x},t)=U\overline{T}-T(DU_{x}+EU_{y}+FU)$ such that $T(t)=\frac{t^3}{6}-t$, and with the order of the derivatives $\alpha_6=\frac{5}{3},\;\alpha_5=\frac{7}{5},\;\alpha_4=\frac{4}{3},\;\alpha_3=\frac{6}{5},\;\alpha_2=\frac{3}{2},\;\alpha_1=\frac{2}{3},$ and $\alpha_0=\frac{1}{2}$, and their coefficients $\gamma_6=\Gamma(\frac{1}{3}),\;\gamma_5=4\pi,\;\gamma_4=2\pi,\;\gamma_3=4\pi,\;\gamma_2=\sqrt{\pi},\;\gamma_1=\frac{\Gamma(\frac{1}{3})}{3},$ $\gamma_0=\frac{\sqrt{\pi}}{2},\;$ $D={\frac{r_{2}^{2}-x^{2}}{r_{1}^{2}-y^{2}}},\;
E=\frac{x^{2}-y^{2}}{r_{1}-r_{2}},$
and $F=\sqrt{x^{2}+y^{2}}+r_{3},$  we can set
\[\begin{gathered}
\overline{T}(t)=t^{\textstyle\frac{5}{3}} \left( \frac{
2\pi}{\Gamma(\textstyle\frac{8}{3})}+\textstyle\frac{\Gamma(\textstyle\frac{1}{3})}{3\Gamma(\textstyle\frac{11}{3})}t\right)
+{t^{\textstyle\frac{3}{2}}}\left(\textstyle\frac{20-8t}{15}\right)
+t\left(\textstyle\frac{9}{4}t^{\frac{1}{3}}+\textstyle\frac{4\pi}{\Gamma(\textstyle\frac{13}{5})}t^{\frac{3}{5}}+\textstyle\frac{4\pi t^{\frac{4}{5}}}{\Gamma(\textstyle\frac{14}{5})}\right) \hfill \\
-\left(t^{\textstyle\frac{-2}{3}}+t^{\textstyle\frac{1}{3}}+t^{\textstyle\frac{1}{2}}+t^{\textstyle\frac{-1}{2}} +\sqrt{3}\Gamma(\textstyle\frac{1}{3})t^{\textstyle\frac{-1}{3}}+\left(\sqrt{10+2\sqrt{5}}\right)\Gamma(\textstyle\frac{2}{5})t^{\textstyle\frac{-2}{5}}+\left(\sqrt{10-2\sqrt{5}}\right)\Gamma(\textstyle\frac{1}{5})t^{\textstyle\frac{-1}{5}} \right),  \hfill
\end{gathered}\]
to find the exact solution as 
$u_{exact}=T(t)U(\textbf{x})$. This problem is solved with $N=216$, $M=741$, and $K=16$, for $t=5$. In Figure~\ref{figEX6}, the distribution of the absolute errors on the domain, $L_{\infty}, MRE,$ and $RMS$ of $u$, $u_{xx}$, and $u_{yy}$ are illustrated.
Figure~\ref{figCondEX6} exhibits the behavior of the condition number matrix $A$ as $r^{-2}\times (K+1)^3$; e.g. when $r=1$, $Cond(\pmb{A})\simeq 0.92 \times (K+1)^3$, and when $K=10$, $Cond(\pmb{A})\simeq 1202.46 \times r^{-2}$.
Among these six examples, an interesting point can be concluded that the error distribution into a plan depends not only on the positions of the nodal points but also on the final time solving the problem; with an increasingly asymmetric discretization, and longer computing time, the error distribution becomes more \say{random.}

\begin{figure}[!ht]
\begin{center}
{\includegraphics[width=1\textwidth]{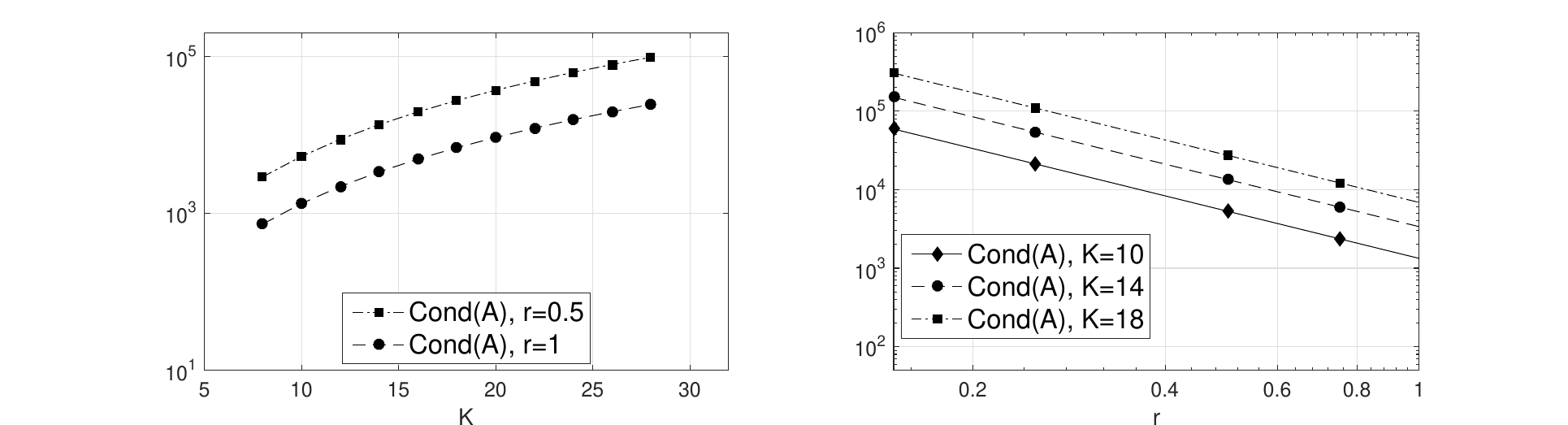}}
\caption{The condition number of matrix $\pmb{A}$ versus the polynomial degree $K$ (left) and versus the mean length of discrete element edges $r$ (right) for Example~\ref{exam6} when $t=5$.}\label{figCondEX6}
\end{center}
\end{figure}

\begin{figure}[!ht]
\begin{center}
\subfigure{\includegraphics[width=.4\textwidth]{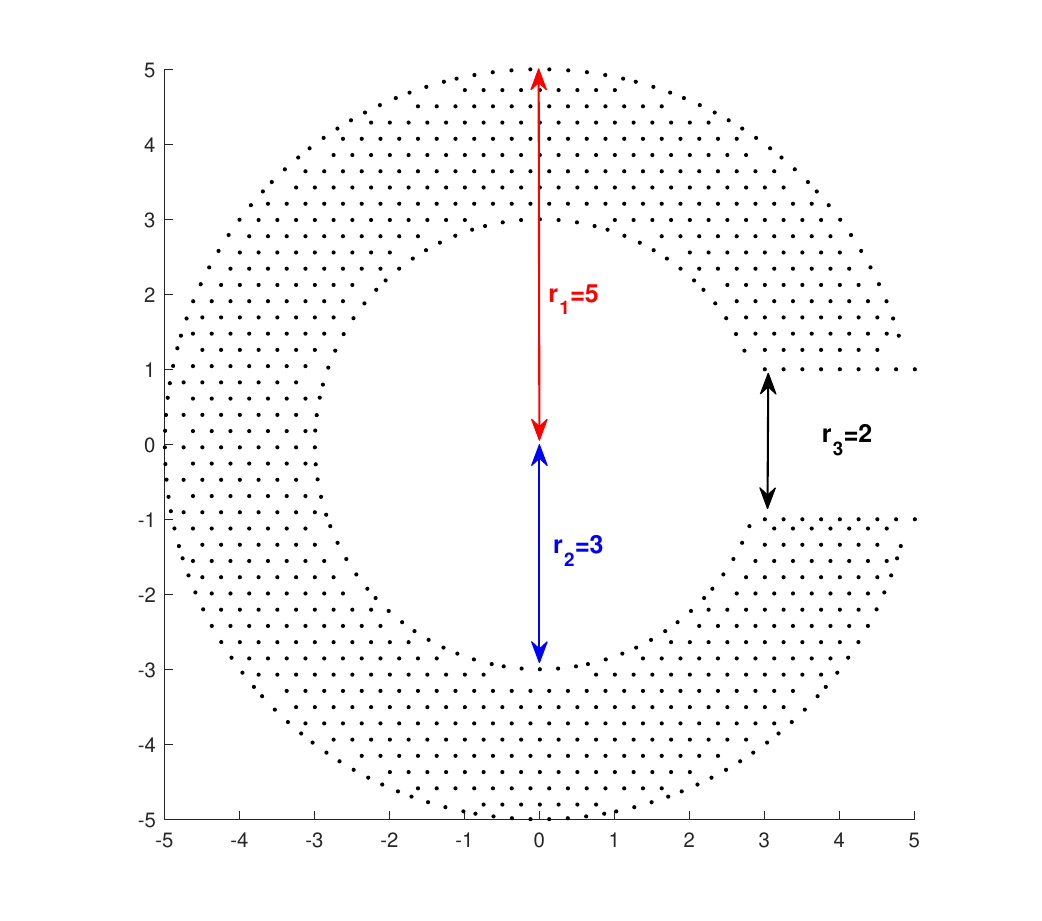}}
\subfigure{\includegraphics[width=.4\textwidth]{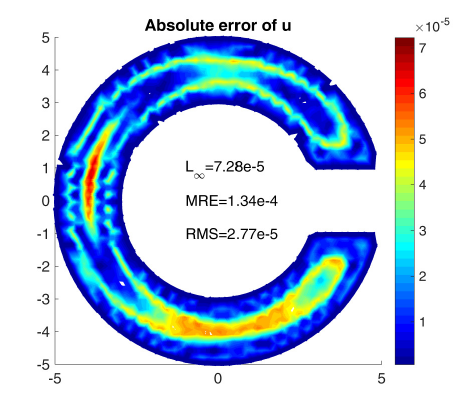}}
\subfigure{\includegraphics[width=.4\textwidth]{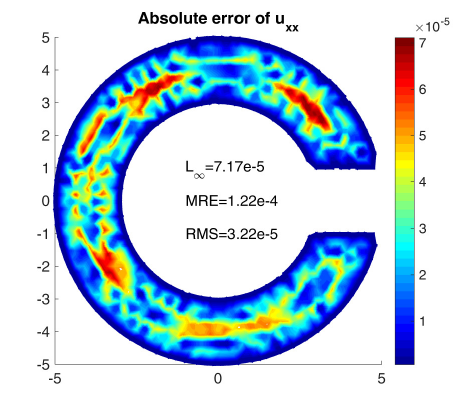}}
\subfigure{\includegraphics[width=.4\textwidth]{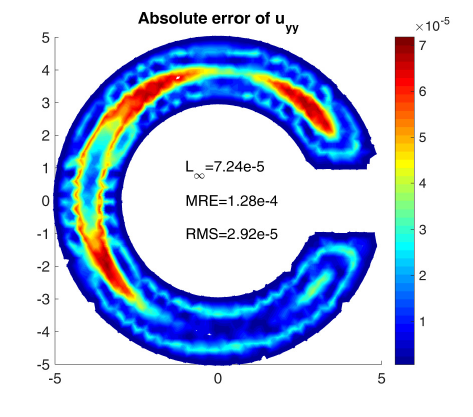}}
\caption{The geometry of the plane body discretization with  $N=216$, $M=741$ with DG-RBF technique and Contour plots of absolute errors when $K=16$ for Example~\ref{exam6}.}\label{figEX6}
\end{center}
\end{figure}

\end{exm}

\section*{Conclusion}
Here, we have proposed a hybrid algorithm to solve two-dimensional multi-order time-fractional partial differential equations. Their general form is given in equations (\ref{eq1}-\ref{eq3}). The method consists of the boundary element method combined with spectral Chebyshev operational matrix. The BEM is used to transfer the corresponding time fractional PDE into a system of ODEs while COM is used to solve the system efficiently. This method is applied to the two-dimensional fractional heat-like, wave-like and diffusion-wave equations, which shows that the errors of the approximate solution decay exponentially. When the exact solution exists, comparison is made with $\left\| {{{\mathbf{U}}_{ex}} - {{\mathbf{u}}_{app}}} \right\|$, and the convergence rate is calculated using Lemma~\ref{lemma1}. When the exact solution is not in hand, the order of convergence is estimated by three approximate solutions with various degrees of Chebyshev polynomials in a same grid-point based on Lemma~\ref{lemma2}. By applying the assumptions of the Lemmas, the numerical results show the efficiency and convergence rate for the proposed hybrid method.
Notwithstanding, it is not easy to emphasize a unique conclusion for the accuracy of the method on the ground that given the vast range of architectures, spectral methods, boundary element methods, fractional calculus, and meshing used with in such a hybrid-technique framework. 
In general, for multi-order two-dimensional time fractional PDE (\ref{eq1}-\ref{eq3}), the current method calculate ${\mathbf{U}}$ for the test examples, with this range of convergence rate:  $2<P_r$-order$<4.5$ for the moderate values of $N$ and $M$. And for multi-term fractional ODE \eqref{eq14} with initial condition \eqref{eq16} the current method works well to calculate solution ${\mathbf{b}}$ with a range of the convergence rate around  $1<P_\tau$-order$<5.5$. Moreover, The condition number of matrix $\pmb{A}$ from linear system~\eqref{matrixform} behaves like $Cond\pmb{A}\simeq{r^{-2}(K+1)^{2}}$ for the problems with $\left\lceil \alpha \right\rceil =1$, and $Cond\pmb{A}\simeq{r^{-2}(K+1)^{3}}$ for the problems with $\left\lceil \alpha \right\rceil =2$. For the future direction, the authors believe that establishing new methods to examine long-term effects of memory in complex systems modeled by fractional calculus is highly required as fractional calculus is a proper mathematical tool for describing memory~\cite{megh}, while the proposed COM technique is not an appropriate scheme for long-term problems. It is instead efficient for problems with multi-term orders. 

\section*{Acknowledgements}
The authors are very grateful to the referees for carefully reading the paper and for their comments and suggestions which have improved the paper.


\newpage
\section*{Appendix}\label{sec: appn}

\subsection*{Algorithm in a nutshell}


\begin{minipage}[t]{0.44\textwidth}
For the convenient, a notation table, and the proposed algorithm's description is given in this section.

At the beginning, determine $N$ boundary points and divide $M$ internal nodal points into groups by the Delaunay graph, then use the RBF method to interpolate the nodal points to its new position. Hence, there is no need to optimize shape parameters for the domain discretization. After discretization, consider $\int_k$ indicating integration on k-element on the boundary (see Figure \ref{fig: discrete} and Table~\ref{tabale: notation}), and set $i, k =1,...,N$, $j=1,...,M$ for the boundary points $x_{BP_i}$, the domain points $x_{IP_j}$, and the algorithm implementation. Notice, different $\theta$ must be considered for computing the boundary integrals at the corner points, particularly in inhomogeneous shapes, in comparing to smooth boundary\cite{katsi2002}.

\end{minipage}\hfill
\begin{minipage}[t]{0.54\textwidth}
  \centering\raisebox{\dimexpr \topskip-\height}{%
  \includegraphics[width=\textwidth]{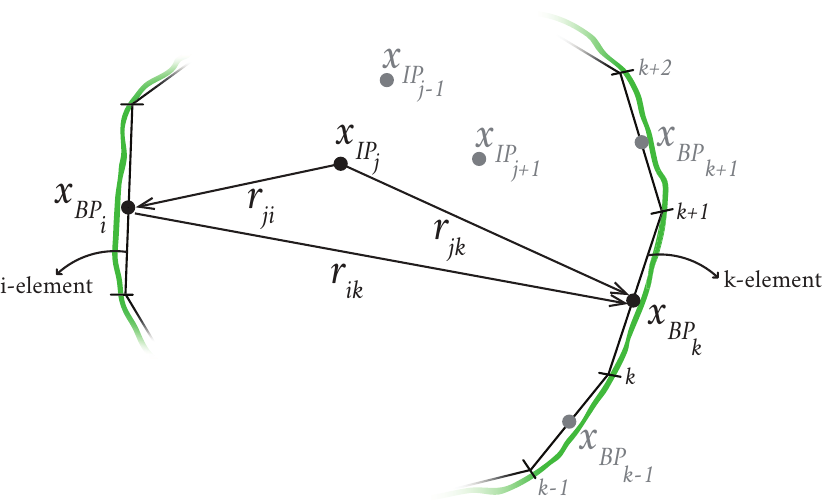}}
  \captionof{figure}{The location of nodal points and relative distances for constant element discretization.}
  \label{fig: discrete}
\end{minipage}

\begin{table}[!ht]
\vspace{.3cm}	
\centering 
\begin{tabular}{l l} 
\noalign{\hrule height 2pt}
 \multicolumn{2}{c}{\textbf{Algorithm}}\\
\noalign{\hrule height 2pt}
\textbf{Step 0:} & {Input $M, N, r_{ik}, r_{jk}, f_{j}(r), \mathbf{\alpha}=\{\alpha_{0},...,\alpha_{k}\}$.}\\
\multicolumn{2}{c}{\textit{Implementation of boundary element method}}\\
\textbf{Step 1:} & {Initialize  $\mathbf{u}^{*}, \mathbf{u}^{*}_{n}, \hat{\mathbf{u}}, \hat{\mathbf{u}}_{n}, \textbf{A}, \textbf{B}, \textbf{C}, \textbf{D}, \textbf{E}, \textbf{F}, \mathbf{\gamma}_{j}, \textbf{h}, \mathbf{\delta}_{1}, \mathbf{\delta}_{2}, \textbf{d}.$}\\ 
 \textbf{Step 2:} & Compute  $\textbf{H}=\tilde{\textbf{H}}-\frac{\theta}{2\pi}\textbf{I}$, where $\textbf{I}$ is an $N \times N$ identity matrix, and $\tilde{\textbf{H}}(i,k)=\int_{k}{\mathbf{u}_n^{*}(r_{ik})ds}$.  \\ 
 \textbf{Step 3:} & Compute ${\textbf{G}}(i,k)=\int_{k}{\mathbf{u}^{*}(r_{ik})ds}$. \\
 \textbf{Step 4:} & Compute  $\bar{\textbf{A}}(i,k)=\frac{1}{2}\hat{\mathbf{u}}(r_{ji})-\displaystyle\sum_{k=1}^{N} \tilde{\textbf{H}}(i,k)\hat {\mathbf{u}} (r_{jk}) + \displaystyle\sum_{k=1} ^ {N} {\textbf{G}} (i,k) \hat {\mathbf{u}}_{n}(r_{jk})$. \\
 \textbf{Step 5:} & Compute $\hat{\textbf{H}}_{pq}(i,k)=\int_{k}{(\mathbf{u}^{*}_{n})_{pq}(r_{ik})ds}$, where $p,\;q=0,\;x,\;y$. \\
 \textbf{Step 6:} & Compute $\hat{\textbf{G}}_{pq}(i,k)=\int_{k}{(\mathbf{u}^{*})_{pq}(r_{ik})ds}$. \\
 \textbf{Step 7:} & Compute $\bar{\textbf{A}}(i,k)=\frac{1}{2}\hat{\mathbf{u}}(r_{ji})-\displaystyle\sum_{k=1}^{N} \tilde{\textbf{H}}(i,k)\hat {\mathbf{u}} (r_{jk}) + \displaystyle\sum_{k=1} ^ {N} {\textbf{G}} (i,k) \hat {\mathbf{u}}_{n}(r_{jk})$. \\
  \textbf{Step 8:} & Construct $\textbf{U}_{pq}$ and $\textbf{c}_{pq}$ by~\eqref{eq:Uc}, consequently, $\textbf{S}$, $\textbf{N}$, and $\textbf{f}(t)$ by~\eqref{eq:SNF}. \\
  \multicolumn{2}{c}{\textit{Solving linear multi-order fractional ODE}}\\
  \textbf{Step 9:} & {Construct vector function $\Phi(t)$ by~\eqref{eq:phi} and \eqref{eq: collocation}, and spectral matrix ${\mathfrak{D}}^{(\upsilon )}$ for $\upsilon\in{\mathbf{\alpha}}$ by~\eqref{eq28}.}\\
   \textbf{Step 10:} & Construct vector function $\textbf{R}(t)$ by~\eqref{eq32}.\\
  \textbf{Step 11} & Construct $M(K+1)$ linear system by solving~\eqref{eq33} and using~\eqref{eq34}. \\
  \textbf{Step 12}  & Solve the generated algebraic system \eqref{matrixform} for the vector $\pmb{\Psi}$. \\
  \textbf{Step 13}  & Compute unknown vector $\textbf{b}(t)$ from~\eqref{eq29}. \\
  \textbf{Step 14}  & Output the corresponding solution of the problem~\eqref{eq1} and its derivatives by applying equation~\eqref{eq12}. \\
  
\hline
\end{tabular}
\label{tabale: algorithm}
\end{table}

\begin{table}[!ht]
\caption{Notation} 
\vspace{.3cm}	
\centering 
\begin{tabular}{p{2cm} p{14cm}} 
\noalign{\hrule height 2pt}
 Symbol & Description \\ 
\noalign{\hrule height 2pt}
 $\Omega$  & Two dimensional domain  \\ 
$\Gamma$  &   Two dimensional boundary \\
$D_c^{\alpha_j}$ & Caputo fractional time derivative of order $\alpha_j$ \\
$u(\textbf{x},t)$ & Unknown field function of spatial $\textbf{x}\in{\Omega\cup\Gamma}$, and time $t$\\
$u_n$  & Normal derivative of $u$  \\
$A, \;B,\; C,$ & \multirow{2}{12cm}{Given coefficient functions of $\textbf{x}$ and $j=0,1,...,k$}\\
$\;D,\; E\;, F,\; \gamma_{j}$ & \\
g  & Given independent function of $\textbf{x}$ and t\\
$\rm B$  & Linear operator with respect to $\textbf{x}$ of order one  \\
$h(\textbf{x},t)$  & Given function in the boundary condition; $\textbf{x}\in{\Gamma}$ \\
$d_{i}(\textbf{x})$  & Given function in the initial condition; $i=0,1,...,m-1$ \\
$\mathfrak{B}(\textbf{x},t)$ &  Unknown fictitious source function \\
$u^*$, $u^{*}_n$  & Fundamental solution of~\eqref{eq5}, and normal derivative of $u^*$ on the boundary\\
$\hat u$, $\hat u_n$ &  Particular solution of~\eqref{eq8}, and normal derivative of $\hat u$\\
 $\varepsilon$  &  Free term coefficient; $\varepsilon=1$ if $\textbf{x}\in{\Omega}$, $\varepsilon=\theta/{2\pi}$ if $\textbf{x}\in{\Gamma}$,  else $\varepsilon=0$, see details in book~\cite{katsi2002}\\
$\theta$ & Interior angle between the tangents of boundary at point $\textbf{x}$, see details in book~\cite{katsi2002}\\
$r$  &  Distance between two points (or mean of all distances of internal points)\\
$M$ & Number of interior points after discretization  \\
${x_{I{P_j}}}$ & $M$ internal nodal points; $j=1,...,M$\\
$N$ & Number of boundary nodal points after discretization \\
${x_{B{P_j}}}$ & $N$ boundary nodal points; $j=1,...,N$\\
$f_{j}(r)$ & Radial basis approximating functions, $j=1,2,...,M$\\
$\mathbf{H}$, $\rm \textbf{G}$ & {$N \times N$ known coefficient matrices  from the integration of the kernel functions on the boundary }\\
${\mathbf{\bar A}}$  & $N \times M$  known coefficient matrix  from the integration of the kernel function on the domain \\
$\textbf{u}, \textbf{u}_n$  & Unknown vectors of the nodal values of the boundary displacements and their normal derivatives \\
$\textbf{b}(t)$ &  Vector of the nodal values of the fictitious source at the $M$ internal nodal points\\
$\mathbf{\delta}_1$, $\mathbf{\delta}_2$ & $N \times N$ known diagonal matrices \\
 $\textbf{h}(t)$ & Known boundary vector   of ${h( {{x_{B{P_j}}},t} )}$, $j=1,...,N$
 \\
${{\bf{\hat u}}_{pq}}$  & Vector of values for $u$ and its derivatives at the $M$ internal nodal points; ${{\bf{\hat u}}_{00}}=\bf{\hat u}$  \\
 ${{\bf{\hat H}}_{pq}}$, ${{\bf{\hat G}}_{pq}}$ &  $M \times N$ known coefficient matrices from the integration of the kernel functions on the boundary  \\
 ${{\bf{\hat A}}_{pq}}$ &  $M \times M$ known coefficient matrix  from the integration of the kernel functions on the domain \\
$\textbf{A}, \;\textbf{B},\; \textbf{C},$ & \multirow{2}{14cm}{$M \times M$ known diagonal matrices including the nodal values of the corresponding functions $A(\bf{x})$, $B(\bf{x})$, $C(\bf{x})$, $D(\bf{x})$, $E(\bf{x})$, $F(\bf{x})$, $\gamma_{j}(\bf{x})$}\\
$\;\textbf{D},\; \textbf{E}\;, \textbf{F},\; \mathbf{\gamma}_{j}$ & \\
$\textbf{g}(t)$ & Known internal vector of  $g(x_{IP_{j}},t),\;j=1,...,M$\\
$T_{L,i}(t)$  & Shifted Chebyshev polynomials of degree $i$ on the interval $t\in{(0,L)}$ \\
${\mathfrak{D}}^{(\upsilon )}$  & $(K+1)\times(K+1)$ Chebychev operational matrix of derivatives of order $\upsilon$ in the sense of Caputo\\
$\pmb{\Psi}$ & $M\times(K+1)$ unknown matrix \\
$\pmb{\omega}$ & $M\times(K+1)$ known matrix\\
$\textbf{R}(t)$  & Residual vector of~\eqref{eq14} with length $M$\\
$L_{\infty}$  & Maximum error  \\
 $MRE$ & Maximum relative error  \\
$RMS$  & Root mean square \\
$P_{r}$-order & convergence order of BEM approximated by the lemma \ref{lemma1} \\
$P_{\tau}$-order & convergence order of COM approximated by the lemma \ref{lemma2} \\
$Cond(\pmb{A})$ & condition number of matrix $\pmb{A}$ of the algebraic system equations \eqref{matrixform}\\
\hline
\end{tabular}
\label{tabale: notation}
\end{table}

\end{document}